\documentclass[12pt]{article}
\usepackage{amsmath,amsthm,amssymb}
\usepackage{graphicx, epsfig}
\usepackage{color}

\makeatletter
\newfont{\footsc}{cmcsc10 at 8truept}
\newfont{\footbf}{cmbx10 at 8truept}
\newfont{\footrm}{cmr10 at 10truept}
\makeatother
\newcommand{\bc}{\begin{center}}
\newcommand{\ec}{\end{center}}
\newcommand{\be}{\begin{equation}}
\newcommand{\ee}{\end{equation}}
\newcommand{\bea}{\begin{eqnarray}}
\newcommand{\eea}{\end{eqnarray}}
\newcommand{\bfl}{\begin{flushleft}}
\newcommand{\efl}{\end{flushleft}}
\newcommand{\bdm}{\begin{displaymath}}
\newcommand{\edm}{\end{displaymath}}
\newcommand{\ba}{\begin{array}}
\newcommand{\ea}{\end{array}}
\newcommand{\bd}{\begin{description}}
\newcommand{\ed}{\end{description}}
\newcommand{\ben}{\begin{enumerate}}
\newcommand{\een}{\end{enumerate}}
\newcommand{\beas}{\begin{eqnarray*}}
\newcommand{\eeas}{\end{eqnarray*}}
\newcommand{\bb}{\begin{thebibliography}{99}}
\newcommand{\eb}{\end{thebibliography}}
\newcommand{\bs}{\begin{sloppypar}}
\newcommand{\es}{\end{sloppypar}}

\newtheorem{thm}{Theorem}[section]
\newtheorem{cor}[thm]{Corollary}
\newtheorem{lem}[thm]{Lemma}
\newtheorem{prop}[thm]{Proposition}

\newtheorem{con}[thm]{Conjecture}  

\theoremstyle{definition}
\newtheorem{defn}[thm]{Definition}

\theoremstyle{remark}

\newcommand{\ncra}[2]{\mbox{$\left[\begin{array}{c}#1\\#2\end{array}\right]$}}

\newcommand{\eval}[2][\right]{\relax
  \ifx#1\right\relax \left.\fi#2#1\rvert}

\title{Kocay's lemma, Whitney's theorem, and some polynomial
invariant reconstruction problems}
\author{Bhalchandra D. Thatte \\
\small Allan Wilson Centre for Molecular Ecology and Evolution, \\[-0.8ex]
\small and Institute of Fundamental Sciences, \\[-0.8ex]
\small Massey University, Palmerston North, New Zealand \\[-0.8ex]
\small \texttt{b.thatte@massey.ac.nz}}

\date{\small 
Submitted: June 29, 2004 }
\begin{document}
\maketitle

\begin{abstract}
Given a graph $G$, an incidence matrix $\mathcal{N}(G)$
is defined on the set of distinct isomorphism
types of induced subgraphs of $G$. It is proved that
Ulam's conjecture is true if and only if the
$\mathcal{N}$-matrix is a complete graph invariant.
Several invariants of a graph are then shown to
be reconstructible from its $\mathcal{N}$-matrix. The
invariants include the characteristic polynomial, the
rank polynomial, the number of spanning trees 
and the number of hamiltonian cycles in a graph.
These results are stronger than the original results of
Tutte in the sense that actual subgraphs are not used.
It is also proved that the characteristic polynomial of
a graph with minimum degree 1 can be computed from the
characteristic polynomials of all its induced proper subgraphs.
The ideas in Kocay's lemma play a crucial role in most proofs.
Kocay's lemma is used to prove Whitney's subgraph
expansion theorem in a simple manner. The reconstructibility
of the characteristic polynomial is then demonstrated
as a direct consequence of Whitney's theorem as formulated here.
\end{abstract}

\section {Introduction}
\label{sec-intro}
Suppose we are given the collection of induced subgraphs
of a graph. There is a natural partial order
on this collection defined by the induced subgraph relationship
between members of the collection.
An incidence matrix may be constructed to represent this
relationship along with the multiplicities
with which members of the collection appear as induced subgraphs
of other members. Given such a matrix, is it possible to construct
the graph or compute some of its invariants? Such a question is
motivated by the treatment of chromatic polynomials
in Biggs \cite{biggs1993}. Biggs demonstrates that it is
possible to compute the chromatic polynomial of a graph from
its incidence matrix. The idea of Kocay's lemma in
graph reconstruction theory is extremely useful
in studying the question for other invariants.
In this paper, we present several results on a relationship between
Ulam's reconstruction conjecture and the incidence matrix.
Extending the reconstruction results of Tutte and Kocay,
we show that many graph invariants can be computed
from the incidence matrix.
We then consider the problem of computing the characteristic polynomial
of a graph from the characteristic polynomials of all induced
proper subgraphs. Finally, we present a new short proof of
Whitney's subgraph expansion theorem, and demonstrate the
reconstructibility of the characteristic polynomial of a graph using
Whitney's theorem.

\subsection {Notation}

We consider only finite simple graphs in this paper.
Let $G$ be a graph with vertex set $VG$ and edge set
$EG$. The number of vertices of $G$ is denoted by $v(G)$
and the number of edges is denoted by $e(G)$.
When $VG=\emptyset $, we denote $G$ by $\Phi$, and call the
graph a {\em null graph}. When $EG=\emptyset $, we call
the graph an {\em empty graph}.
When $F$ is a subgraph of $G$, we write $F\subseteq G$,
and when $F$ is a proper subgraph of $G$, we write $F\subsetneq G$.
The subgraph of $G$ induced by $S \subseteq VG$
is the subgraph whose vertex set is $S$ and whose edge
set contains all the edges having both end vertices
in $S$. It is denoted by $G_S$. The subgraph of $G$ induced by
$VG - S$ is denoted by  $G-S$, or simply $G-u$ if $S=\{u\}$. 
A subgraph of $G$ with vertex set $V\subseteq VG$ and edge set
$E\subseteq EG$ is denoted by $G_{(V,E)}$, or just $G_E$
if $V$ consists of the end vertices of edges in $E$.
The same notation is used when $E = (e_1,e_2, \ldots ,e_k)$
is a tuple of edges, some of which may be identical.
Isomorphism of two graphs $G$ and $H$ is
denoted by $G\cong H$. For $i > 0$, a graph isomorphic to a
cycle of length $i$ is denoted by $C_i$,
and the number of cycles of length $i$ in $G$ is
denoted by $\psi_i(G)$, where, as a convention,
$C_i \cong K_i$ for $i \in \{1,2\}$. The number of hamiltonian
cycles is denoted by a special symbol $ham(G)$ instead of
$\psi_{v(G)}(G)$.
While counting the number of subgraphs of a graph $G$
that are isomorphic to a graph $F$, it is important to make
a distinction between induced subgraphs and edge subgraphs.
The number of subgraphs of $G$ that are isomorphic to $F$
is denoted by $\ncra{G}{F}$, and the number of induced
subgraphs of $G$ that are isomorphic to $F$ is denoted by
$\displaystyle\binom{G}{F}$. The two numbers are related by
\begin{equation}
\ncra{G}{F} = \sum_{H|VH=VF}\displaystyle\binom{G}{H}\ncra{H}{F}
\end{equation}
where the summation is over distinct isomorphism types
of graphs $H$.
The characteristic polynomial of $G$ is denoted by
$P(G;\lambda) = \sum_{i = 0}^{v(G)} c_i(G) \lambda^{v(G)-i}$.
The collection 
$\mathcal{PD}(G) = \{P(G-S;\lambda)\mid S\subsetneq VG\}$
is called the {\em complete polynomial deck} of $G$.
Note that a polynomial may appear in the collection more than once.
The {\em rank} of a graph $G$, which has $comp(G)$ components,
is defined by $v(G) - comp(G)$, and its {\em co-rank} is defined
by $e(G)-v(G) + comp(G)$. The rank polynomial of $G$
is defined by $R(G;x,y) = \sum \rho_{rs}x^ry^s$,
where $\rho_{rs}$ is the number of subgraphs of $G$
with rank $r$ and co-rank $s$. The set of consecutive
integers from $a$ to $b$ is denoted by $[a,b]$;
in particular, $N_k = [1,k]$.

\subsection {Ulam's Conjecture}

The {\em vertex deck} of a graph $G$ is the collection
$\mathcal{VD}(G) = \{G-v\mid v \in VG\}$, where the subgraphs
in the collection are `unlabelled' (or isomorphism types).
Note that the vertex deck is not exactly a set: an
isomorphism type may appear more than once in the vertex deck.
A Graph $G$ is said to be {\em reconstructible} if its
isomorphism class is determined by $\mathcal{VD}(G)$.
Ulam \cite{ulam1960} proposed the following conjecture.
\begin{con} 
Graphs on more than 2 vertices are reconstructible.
\end{con}

A property or an invariant of a graph $G$ is said to be
reconstructible if it can be calculated from $\mathcal{VD}(G)$.
For example, Kelly's Lemma allows us to count the number of
vertex-proper subgraphs of $G$ of any given type.

\begin{lem} \label{lem-kelly}{\bf (Kelly's Lemma \cite{kelly1957})}
If $F$ is a graph such that $v(F) < v(G)$ then
\begin{equation}
\label{eq-kelly}
\ncra{G}{F} = \frac{1}{v(G)-v(F)}\sum_{u\in VG}\ncra{G-u}{F}
\end{equation}
therefore, $\ncra{G}{F}$ is reconstructible from $\mathcal{VD}(G)$.
Also, in Equation~(\ref{eq-kelly}), 
$\ncra{G}{F}$ and $\ncra{G-u}{F}$ may be replaced by 
$\displaystyle\binom{G}{F}$ and $\displaystyle\binom{G-u}{F}$, respectively.
\end{lem}

Tutte \cite{tutte1977}, \cite{tutte1984} proved the
reconstructibility of the characteristic polynomial
and the chromatic polynomial.
Tutte's results were simplified by an elegant counting argument
by Kocay \cite{kocay1981}. This argument is useful to count
certain subgraphs that span $VG$.

Let $S = \{F_1, F_2, \ldots ,F_k\}$ be a family of graphs.
Let $c(S,H)$ be the number of tuples 
$(X_1, X_2, \ldots ,X_k)$ of subgraphs of $H$ such that
$X_i\cong F_i\, \forall \, i$, and $\cup_{i = 1}^{k}X_i = H$.
We call it the number of {\em $S$-covers} of $H$.

\begin{lem} \label{lem-kocay} {\bf (Kocay's Lemma \cite{kocay1981})}
\begin{equation}
\prod_{i = 1}^{k} \ncra{G}{F_i} = \sum_{X} c(S,X)\ncra{G}{X}
\end{equation}
where the summation is over all isomorphism
types of subgraphs of $G$.
Also,  if $v(F_i)\, <\, v(G)\,\forall\, i$ then
$ \sum_{X} c(S,X)\ncra{G}{X}$ over all isomorphism types $X$ of
spanning subgraphs of $G$ can be reconstructed
from the vertex deck of $G$.
\end{lem}

We refer to \cite{bondy1991} for a survey of reconstruction problems.

\subsection {The chromatic polynomial and the $\mathcal{N}$-matrix}

Stronger reconstruction results on the chromatic polynomial were implicit
in Whitney's work \cite{whitney1932}, although Ulam's
conjecture had not been posed at the time.
Motivation for some of the work presented in this paper comes from
Whitney's work on the chromatic polynomials.
The discussion of the chromatic polynomial presented here
is based on \cite{biggs1993}.

A graph $G$ is called {\em quasi-separable}
if there exists $K\subsetneq VG$ such that $G_K$ is a complete graph
and $G-K$ is disconnected. If $|K| \leq 1$ then $G$ is
said to be separable.
\begin{thm} (Theorem 12.5 in \cite{biggs1993})
\label{thm-biggs}
The chromatic polynomial of a graph is determined by its proper
induced subgraphs that are not quasi-separable.
\end{thm}

The procedure of computing the chromatic polynomial may be
outlined as follows.
First a matrix  $\mathcal{N}(G) = (N_{ij})$ is constructed.
The rows and the columns of $\mathcal{N}(G)$ are indexed by
induced subgraphs $\Lambda_1, \Lambda_2, \ldots , \Lambda_I = G$,
which are the distinct isomorphism types of 
non-quasi-separable induced subgraphs of $G$.
The list includes $K_1 = \Lambda_1$ and $K_2=\Lambda_2$.
The indexing graphs are ordered in such a way that $v(\Lambda_i)$
are in non-decreasing order. The entry $N_{ij}$ is the number of
induced subgraphs of $\Lambda_i$ that are isomorphic to $\Lambda_j$.
It is a lower triangular matrix with diagonal entries 1.
The computation of the chromatic polynomial is performed
by a recursive procedure beginning with the first row
of the $\mathcal{N}$-matrix, computing at each step
certain polynomials in terms of the corresponding polynomials
for non-quasi-separable induced subgraphs on fewer vertices.
A few observations about the procedure are useful 
to motivate the work in this paper.
The graphs $C_4$ and $K_4$ are the only non-quasi-separable
graphs on $4$ vertices.
Also, for any $i$, $N_{i1}=v(\Lambda_i)$,
and $N_{i2}= e(\Lambda_i)$. Therefore, graphs
on 4 or fewer vertices that index the first few rows
of the $\mathcal{N}$-matrix can be inferred from the matrix entries.
Therefore, we conclude that the computation of the chromatic
polynomial can be performed on the matrix entries alone, even if the induced
subgraphs indexing the rows and the columns of $\mathcal{N}(G)$
are unspecified. Therefore, we will think of
the $\mathcal{N}$-matrix as {\em unlabelled }, that is, we
will assume that the induced subgraphs indexing the rows and the
columns are not given.

A natural question is what other invariants
can be computed from the (unlabelled) $\mathcal{N}$-matrix? 
Obviously, the characteristic polynomial $P(G;\lambda)$
cannot always be computed from $\mathcal{N}(G)$.
For example, the only non-quasi-separable
induced subgraphs of any tree $T$ are $K_1$ and $K_2$,
so $P(T;\lambda)$ cannot be computed from $\mathcal{N}(T)$.
Therefore, we omit the restriction of
non-quasi-separability on the induced subgraphs used in the
construction of the incidence matrix. We then investigate which
invariants of a graph are determined by its $\mathcal{N}$-matrix.

The Sections ~\ref{sec-n-rec1} and ~\ref{sec-nmatrix}
are devoted to the study of reconstruction from the
$\mathcal{N}$-matrix.
In Section~\ref{sec-n-rec1}, we formally define the
$\mathcal{N}$-matrix, and the related concept of
the edge \text{labelled} poset of induced subgraphs of a graph.
We then prove several basic results on the relationship
between the $\mathcal{N}$-matrix, the edge labelled poset
and reconstruction. In particular, we show that
Ulam's conjecture is true if and only if the 
$\mathcal{N}$-matrix itself is a complete graph invariant.
We then prove that Ulam's conjecture is true if and only if
the edge labelled poset has no non-trivial automorphisms.
We also prove the $\mathcal{N}$-matrix reconstructibility of
trees and forests.

In Section~\ref{sec-nmatrix} we compute several invariants of a graph
from its $\mathcal{N}$-matrix.
We prove that the characteristic polynomial $P(G;\lambda)$ of
a graph $G$, its rank polynomial $R(G;x,y)$,
the number of spanning trees in $G$, the number of
Hamiltonian cycles in $G$ etc., can be computed from $\mathcal{N}(G)$.
In the standard proof of the reconstructibility of these
invariants, one first counts the disconnected subgraphs of
each type, (see \cite{bondy1991}). In view of
Theorem~\ref{thm-equiv-disc}, the proofs in Section~\ref{sec-nmatrix}
are more involved. Theorem~\ref{thm-equiv-disc} implies that
if there are counter examples to Ulam's conjecture then
there are many more counter examples to reconstruction
from the $\mathcal{N}$-matrix. Therefore, we
hope that the study of $\mathcal{N}$-matrix reconstructibility
will highlight new difficulties. Similar generalisations
of the reconstruction problem were also suggested by
Tutte, (notes on pp. 123-124 in \cite{tutte1984}).

\subsection {Reconstruction of the characteristic polynomial}

The proof of the reconstructibility of the characteristic
polynomial of a graph from its $\mathcal{N}$-matrix
is also of independent technical interest, since other
authors have considered the question of computing $P(G; \lambda)$ given
the {\em polynomial deck} $\{P(G-u;\lambda); u \in VG\}$.
This question was originally proposed by Gutman and
Cvetkovic \cite{gut-cve1975}, and has been studied by
others, for example, \cite{schwenk1979} \& \cite{sciriha2002}.
This question remains open.
So we consider a weaker question in Section~\ref{sec-poly}:
the question of computing the characteristic polynomial of a graph
from its complete polynomial deck. Here we present basic facts about
the characteristic polynomial, and outline the idea of
Section~\ref{sec-poly}.

\begin{defn}
A graph is called {\em elementary} if each of its
components is $1$-regular or $2$-regular.
In other words, each component of an elementary graph is
a single edge ($K_2$) or a cycle ($C_r; r >2$).
\end{defn}

Let $L_i$ be the collection of all unlabelled $i$-vertex
elementary graphs. So, $L_0 = \{\Phi \}$, $L_1 = \emptyset$,
$L_2 = \{K_2\}$, and so on.

\begin{lem}
\label{lem-sachs} (Proposition 7.3 in \cite{biggs1993})
Coefficients of the characteristic polynomial of
a graph $G$ are given by
\begin{equation}
(-1)^i c_i(G) = \sum_{F \in L_i, F\subseteq G} (-1)^{r(F)} 2^{s(F)}
\end{equation}
where $r(F)$ and $s(F)$ are the rank and the co-rank of
$F$, respectively.
\end{lem}

Thus, $c_0(G) = 1$, $c_1(G) = 0$, and $c_2(G) = e(G)$.

\begin{lem} (Note 2d in \cite{biggs1993})
\label{lem-derivative}
Let $P'(G;\lambda)$ denote the first derivative
of $P(G;\lambda)$ with respect to $\lambda$. Then,
\begin{equation}
\label{eq-derivative}
P'(G;\lambda) = \sum_{u\in VG}P(G-u;\lambda)
\end{equation}
\end{lem}

From the above two lemmas, it is clear that the problem of
reconstructing a characteristic polynomial (either from the
vertex deck or the complete polynomial deck) reduces to
computing the coefficient $c_{v(G)}(G)$,
which is the constant term in $P(G;\lambda)$.
This in turn is a problem of counting the
elementary spanning subgraphs of $G$ - a problem that
can be solved using Kocay's Lemma in case of reconstruction
from the vertex deck. Motivated by Kocay's Lemma, we ask the
following question. Suppose the coefficients 
$c_{i_1}(G), c_{i_1}(G), \ldots , c_{i_k}(G)$ are known,
and $i_1 + i_2 + \ldots + i_k \geq v(G)$. If the coefficients
$c_{i_j}; 1 \leq j \leq k$ are multiplied, can we get some information
about the spanning subgraphs of $G$? This is especially tempting
if $i_1+i_2 + \ldots + i_k = v(G)$, since the product is expected
to have some relationship with the disconnected spanning
elementary subgraphs of $G$. This idea is explored in
Section~\ref{sec-poly}.

In Section~\ref{sec-whitney}, we present a very simple new proof of
Whitney's subgraph expansion theorem, again based on Kocay's
lemma. We then present a more direct argument to compute 
the characteristic polynomial of a graph from its vertex deck,
based on our formulation of Whitney's theorem.

\section {Ulam's conjecture and the $\mathcal{N}$-matrix}
\label{sec-n-rec1}


Let $\Lambda(G) = \{\Lambda_i; \,i \in [1,I]\}$ be the set of distinct
isomorphism types of nonempty induced subgraphs of $G$.
We call this the $\Lambda$-deck of $G$.
Let $\mathcal{N}(G)=(N_{ij})$ be an $I$ x $I$ incidence matrix
where $N_{ij}$ is the number of induced subgraphs
of $\Lambda_i$ that are isomorphic to $\Lambda_j$. Thus
$N_{ii}$ is 1 for all $i \in [1,I]$.
We call an invariant of a graph 
{\em $\mathcal{N}$-matrix reconstructible}
if it can be computed from the (unlabelled) $\mathcal{N}$-matrix of the graph.


As an example, the ladder graph $L_3$ and its collection
of distinct induced subgraphs with nonempty edge sets
are shown in Figure 1. Below each graph (except $L_3$) is shown its
multiplicity as an induced subgraph in $L_3$ and its name.

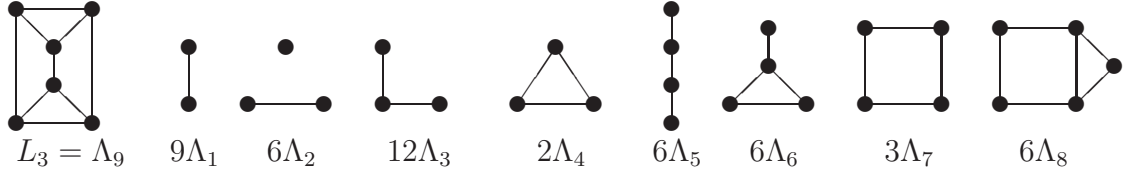
\begin{figure}[ht]
\label{fig-l3}
\vspace*{0.2in}
\setlength{\unitlength}{0.0505mm}
\noindent
\begin{picture}(200,400)
\put(0,100){\circle*{40}}
\put(200,100){\circle*{40}}
\put(0,400){\circle*{40}}
\put(200,400){\circle*{40}}
\put(100,200){\circle*{40}}
\put(100,300){\circle*{40}}
\put(0,100){\line(1,0){200}}
\put(0,100){\line(0,1){300}}
\put(0,400){\line(1,0){200}}
\put(200,100){\line(0,1){300}}
\put(0,100){\line(1,1){100}}
\put(200,100){\line(-1,1){100}}
\put(0,400){\line(1,-1){100}}
\put(200,400){\line(-1,-1){100}}
\put(100,200){\line(0,1){100}}
\put(0,0){$L_3 = \Lambda_9$}
\end{picture}
\hspace*{0.2in}
\begin{picture}(100,400)
\put(100,150){\circle*{40}}
\put(100,300){\circle*{40}}
\put(100,150){\line(0,1){150}}
\put(50,0){9$\Lambda_1$}
\end{picture}
\hspace*{0.2in}
\begin{picture}(200,400)
\put(0,150){\circle*{40}}
\put(200,150){\circle*{40}}
\put(100,300){\circle*{40}}
\put(0,150){\line(1,0){200}}
\put(50,0){6$\Lambda_2$}
\end{picture}
\hspace*{0.2in}
\begin{picture}(200,400)
\put(0,150){\circle*{40}}
\put(150,150){\circle*{40}}
\put(0,300){\circle*{40}}
\put(0,150){\line(1,0){150}}
\put(0,150){\line(0,1){150}}
\put(10,0){12$\Lambda_3$}
\end{picture}
\hspace*{0.2in}
\begin{picture}(200,400)
\put(0,150){\circle*{40}}
\put(200,150){\circle*{40}}
\put(100,300){\circle*{40}}
\put(0,150){\line(1,0){200}}
\put(0,150){\line(2,3){100}}
\put(200,150){\line(-2,3){100}}
\put(50,0){2$\Lambda_4$}
\end{picture}
\hspace*{0.2in}
\begin{picture}(50,400)
\put(50,100){\circle*{40}}
\put(50,200){\circle*{40}}
\put(50,300){\circle*{40}}
\put(50,400){\circle*{40}}
\put(50,100){\line(0,1){100}}
\put(50,200){\line(0,1){100}}
\put(50,300){\line(0,1){100}}
\put(0,0){6$\Lambda_5$}
\end{picture}
\hspace*{0.2in}
\begin{picture}(200,400)
\put(0,150){\circle*{40}}
\put(200,150){\circle*{40}}
\put(100,250){\circle*{40}}
\put(100,350){\circle*{40}}
\put(0,150){\line(1,0){200}}
\put(0,150){\line(1,1){100}}
\put(200,150){\line(-1,1){100}}
\put(100,250){\line(0,1){100}}
\put(50,0){6$\Lambda_6$}
\end{picture}
\hspace*{0.2in}
\begin{picture}(200,400)
\put(0,150){\circle*{40}}
\put(200,150){\circle*{40}}
\put(0,350){\circle*{40}}
\put(200,350){\circle*{40}}
\put(0,150){\line(1,0){200}}
\put(0,150){\line(0,1){200}}
\put(0,350){\line(1,0){200}}
\put(200,150){\line(0,1){200}}
\put(50,0){3$\Lambda_7$}
\end{picture}
\hspace*{0.2in}
\begin{picture}(200,400)
\put(0,150){\circle*{40}}
\put(200,150){\circle*{40}}
\put(0,350){\circle*{40}}
\put(200,350){\circle*{40}}
\put(300,250){\circle*{40}}
\put(0,150){\line(1,0){200}}
\put(0,150){\line(0,1){200}}
\put(0,350){\line(1,0){200}}
\put(200,150){\line(0,1){200}}
\put(200,150){\line(1,1){100}}
\put(200,350){\line(1,-1){100}}
\put(50,0){6$\Lambda_8$}
\end{picture}
\caption[]{$L_3$ and its induced subgraphs.}
\end{figure}

\noindent The rows and the columns of $\mathcal{N}$-matrix of $L_3$
are both indexed by $\Lambda_1$ to $\Lambda_9$.
The $\mathcal{N}$-matrix of $L_3$ is shown below.
\begin{equation}
\mathcal{N}(L_3) = 
\begin{pmatrix}
1 & 0 & 0 & 0 & 0 & 0 & 0 & 0 & 0 \\
1 & 1 & 0 & 0 & 0 & 0 & 0 & 0 & 0 \\
2 & 0 & 1 & 0 & 0 & 0 & 0 & 0 & 0 \\
3 & 0 & 0 & 1 & 0 & 0 & 0 & 0 & 0 \\
3 & 2 & 2 & 0 & 1 & 0 & 0 & 0 & 0 \\
4 & 1 & 2 & 1 & 0 & 1 & 0 & 0 & 0 \\
4 & 0 & 4 & 0 & 0 & 0 & 1 & 0 & 0 \\
6 & 3 & 6 & 1 & 2 & 2 & 1 & 1 & 0 \\
9 & 6 & 12 & 2 & 6 & 6 & 3 & 6 & 1 \\
\end{pmatrix}
\end{equation}


Let us associate an {\em edge labelled poset}
with the graph $G$. Define a partial order $\preceq $
on the set $\Lambda(G)$ as follows:
$\Lambda_j \preceq \Lambda_k$ if and only if $\Lambda_j$ is an
induced subgraph of $\Lambda_k$. 
This poset is denoted by $(\Lambda(G), \preceq)$.
We make the poset $(\Lambda(G), \preceq)$
an edge labelled poset by assigning a positive integer to
every edge of its Hasse diagram, such that if $\Lambda_k$
covers $\Lambda_j$ then the edge label on $\Lambda_j$-$\Lambda_k$ is 
$\displaystyle\binom{\Lambda_k}{\Lambda_j}$.
We say that two edge labelled posets are isomorphic
if they are isomorphic as posets, and there is an isomorphism
between them that preserves the edge labels. This naturally
leads to the notion of the {\em abstract edge labelled poset}
of $G$: it is the isomorphism class of the edge labelled poset
of $G$. Note that the notion of the abstract edge labelled
poset of a graph is not to be confused with the isomorphism
class of the Hasse diagram as a graph.
An isomorphism from an edge labelled poset to itself 
is called an automorphism of the edge labelled poset.
We denote the abstract edge labelled poset of $G$ by
$\mathcal{ELP}(G)$. The Hasse diagram of the abstract
edge labelled poset is simply the Hasse diagram of the
edge labelled poset of $G$ with labels $\Lambda_i$ removed.
The Hasse diagram of $\mathcal{ELP}(L_3)$ is shown in Figure 2.

\vspace*{0.2in}
\begin{figure}[ht]
\label{fig-hasse}
\setlength{\unitlength}{0.0505mm}
\noindent
\begin{center}
\begin{picture}(600,1300)
\put(300,100){\circle*{40}}
\put(0,400){\circle*{40}}
\put(300,400){\circle*{40}}
\put(600,400){\circle*{40}}
\put(0,700){\circle*{40}}
\put(300,700){\circle*{40}}
\put(600,700){\circle*{40}}
\put(300,1000){\circle*{40}}
\put(300,1300){\circle*{40}}
\put(300,100){\line(-1,1){300}}
\put(100,200){2}
\put(300,100){\line(0,1){300}}
\put(470,200){1}
\put(300,100){\line(1,1){300}}
\put(320,250){3}
\put(0,400){\line(0,1){300}}
\put(-80,500){4}
\put(0,400){\line(1,1){300}}
\put(100,550){2}
\put(0,400){\line(2,1){600}}
\put(120,400){2}
\put(300,400){\line(0,1){300}}
\put(320,450){1}
\put(600,400){\line(0,1){300}}
\put(620,500){2}
\put(600,400){\line(-1,1){300}}
\put(450,465){1}
\put(0,700){\line(1,1){300}}
\put(50,830){1}
\put(300,700){\line(0,1){300}}
\put(320,800){2}
\put(600,700){\line(-1,1){300}}
\put(530,830){2}
\put(300,1000){\line(0,1){300}}
\put(320,1100){6}
\end{picture}
\end{center}
\caption[]{The abstract edge labelled poset of $L_3$.}
\end{figure}
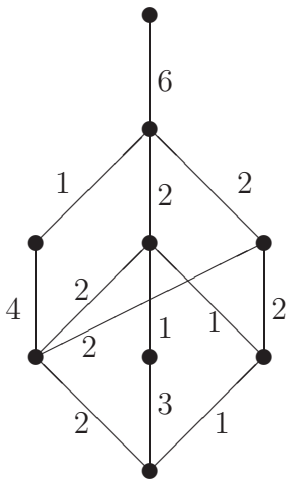

\begin{lem}
\label{lem-rank}
There is a rank function on $\rho $ on $\mathcal{ELP}(G)$ such that
$\rho(\Lambda_i) = \rho(\Lambda_j) +1$ whenever
$\Lambda_i $ covers $\Lambda_j$.
\end{lem}
\begin{proof}
Each $\Lambda_i$ in $\Lambda(G)$ is nonempty. Therefore,
for each $\Lambda_i$ in $\Lambda(G)$ and for each
$k$ such that $2 \leq k \leq v(\Lambda_i)$ there
is at least one nonempty induced subgraph $\Lambda_j$ of
$\Lambda_i$ such that $v(\Lambda_j) = k$. Moreover,
empty induced subgraphs do not belong to $\Lambda(G)$. Therefore,
$\rho(\Lambda_i) = v(\Lambda_i)$ meets the requirements
of a rank function.
\end{proof}

Stanley \cite{stanley1997} defines a rank function such that
the $\rho(x) = 0$ for a minimal element $x$. But we have
deviated from that convention since $\rho(\Lambda_i) = v(\Lambda_i)$
for each $\Lambda_i \in \Lambda(G)$ is more convenient here.
We now demonstrate that $\mathcal{N}(G)$ and $\mathcal{ELP}(G)$
are really equivalent, that is, they can be constructed from
each other. 

\begin{lem}
\label{lem-kelly1} Let $F$ and $H$ be two graphs,
and let $q$ be an integer such that $v(F) \leq q \leq v(H)$. Then
\begin{equation}
\label{eq-kelly1}
\sum_{X| v(X) = q}\displaystyle\binom{H}{X}\displaystyle\binom{X}{F}
= \displaystyle\binom{v(H)-v(F)}{q-v(F)}\displaystyle\binom{H}{F} 
\end{equation}
where the summation is over distinct isomorphism types $X$.
\end{lem}

\begin{proof} This is similar to Kelly's Lemma ~\ref{lem-kelly}.
Each induced subgraph of $H$ that is isomorphic to $F$
is also an induced subgraph of $\displaystyle\binom{v(H)-v(F)}{q-v(F)}$ induced
subgraphs of $H$ that have $q$ vertices.
\end{proof}


\begin{lem}
\label{lem-elp_eq_n}
The structures $\mathcal{N}(G)$ and $\mathcal{ELP}(G)$
can be constructed from each other.
\end{lem}
\begin{proof}
We first show how $\mathcal{N}(G)$ is constructed
from $\mathcal{ELP}(G)$. The matrix $\mathcal{N}(G)$ is
an $I\times I$ matrix where $I$ is the number of points
in $\mathcal{ELP}(G)$.
Without the loss of generality, suppose that 
the points of $\mathcal{ELP}(G)$ are labelled from
$\Lambda_1$ to $\Lambda_I$ such that
if $\rho(\Lambda_i) < \rho(\Lambda_j)$ then $i < j$,
where $\rho$ is the rank function defined in Lemma~\ref{lem-rank}.
Correspondingly, the rows and the columns of $\mathcal{N}(G)$
are indexed from $\Lambda_1$ to $\Lambda_I$.
The edge labels in $\mathcal{ELP}(G)$
immediately give some of the entries in $\mathcal{N}(G)$:
if $\Lambda_i $ covers $\Lambda_j$ then 
$N_{ij}$ is the label on the edge joining $\Lambda_i $ and
$\Lambda_j$. The diagonal entries are 1.
Except $N_{11}$, all the other entries in the first row are 0.
We construct the remaining entries of $\mathcal{N}(G)$ by induction
on the rank. The base case is rank 2. It corresponds to
the first row, and is already filled. 
Let $f(r)$ denote the number of points of $\mathcal{ELP}(G)$
that have rank at most $r$.
Suppose now that the first $f(r)$ rows of $\mathcal{N}(G)$
are filled for some $r \geq 2$. Let $\Lambda_i$ be a
graph of rank $r+1$, and let $\Lambda_j$ be a graph of
rank at most $r$. Then $N_{ij}$ is computed by 
applying Lemma~\ref{lem-kelly1} with $q = r$.
\begin{equation}
\sum_{\Lambda_k| \rho(\Lambda_k) = r}
\displaystyle\binom{\Lambda_i}{\Lambda_k} 
\displaystyle\binom{\Lambda_k}{\Lambda_j}
= \displaystyle\binom{v(\Lambda_i)-v(\Lambda_j)}{r-v(\Lambda_j)}
\displaystyle\binom{\Lambda_i}{\Lambda_j}
= (r+1-v(\Lambda_j))N_{ij}
\end{equation}
On the LHS, $\displaystyle\binom{\Lambda_k}{\Lambda_j}$ are known by
induction hypothesis.
Since $\Lambda_k$ are the graphs covered by $\Lambda_i$,
$\displaystyle\binom{\Lambda_i}{\Lambda_k}$ are the edge labels.
Therefore, $N_{ij}$ can be computed.
This completes the construction of $\mathcal{N}(G)$ from $\mathcal{ELP}(G)$.

To construct $\mathcal{ELP}(G)$ from $\mathcal{N}(G)$,
define a partial order $\preceq $ on
$\{\Lambda_1, \Lambda_2, \ldots , \Lambda_I\}$
as follows: $\Lambda_j \preceq \Lambda_i$ if
$N_{ij} \neq 0$. In this poset, if $\Lambda_i$ covers
$\Lambda_j$ then assign an edge label $N_{ij}$
to the edge $\Lambda_j-\Lambda_i$ of the Hasse
diagram of the poset. This completes the construction of
$\mathcal{ELP}(G)$ from $\mathcal{N}(G)$. 

\end{proof}

\begin{lem}
\label{lem-num-vertices}
Given $\mathcal{N}(G)$, $v(\Lambda_i)$ and $e(\Lambda_i)$ can
be counted for each graph in $\Lambda(G)$.
\end{lem}
\begin{proof}
There is a unique row in $\mathcal{N}(G)$
that has only one nonzero entry (the diagonal entry 1).
This row corresponds to $\Lambda_1 \cong K_2$, and
we assume it to be the first row. Now $e(\Lambda_i) = N_{i1}$
for each $\Lambda_i$.

By Lemma~\ref{lem-elp_eq_n}, $\mathcal{ELP}(G)$ is uniquely
constructed. By Lemma~\ref{lem-rank}, the rank function of
the poset defined by $\rho(\Lambda_1)=2$
gives $v(\Lambda_i)=\rho(\Lambda_i)$ for each $\Lambda_i$.
\end{proof}

Now on, without the loss of generality, we will assume that
the nonisomorphic induced subgraphs $\Lambda_1, \Lambda_2, \ldots,
\Lambda_I$ of a graph $G$ under consideration are ordered 
so that $v(\Lambda_i)$ are in a non-decreasing order.
The first row will correspond to $\Lambda_1 \cong K_2$
and the last row to $\Lambda_I \cong G$.


\begin{lem}
\label{lem-kelly3}
The collection $\{\mathcal{N}(G-u)| u \in VG, e(G-u) > 0\}$
is unambiguously determined by $\mathcal{N}(G)$.

\end{lem}

Note that this collection is a ``multiset'', that is,
an $\mathcal{N}$-matrix may appear multiple times in the
collection.

\begin{proof}
Let $j \neq I$. The graph $\Lambda_j$ is a vertex deleted subgraph
of $G$ if and only if for all $i \neq j \neq I$,
$N_{ij}=0$. Now $N(\Lambda_j)$ is obtained by
deleting $k$'th row and $k$'th column for each $k$
such that $N_{jk} = 0$. A multiplicity $N_{Ij}$ is
assigned to $N(\Lambda_j)$.
Equivalently, we can construct $\mathcal{ELP}(G)$
by Lemma~\ref{lem-elp_eq_n}, then construct the down set
$\mathcal{ELP}(\Lambda_j)$ of each $\Lambda_j$ that
is covered by $\Lambda_I = G$, and then 
construct $\mathcal{N}(\Lambda_j)$, and assign it
a multiplicity equal to the edge label on $\Lambda_I-\Lambda_j$. 
\end{proof}

\noindent {\bf Remark} It is is possible that for distinct
$j$ and $k$, the matrices $N(\Lambda_j)$ and $N(\Lambda_k)$
are equal. In this case a multiplicity $N_{Ij}$ is
assigned to $N(\Lambda_j)$ and $N_{Ik}$ is
assigned to $N(\Lambda_k)$ while constructing the above
collection.

\begin{lem}
\label{lem-empty}
Let $rK_1$ be the $r$-vertex empty graph. The number of induced
subgraphs of $G$ isomorphic to $rK_1$ is determined by
$\mathcal{N}(G)$.
\end{lem}
\begin{proof} The required number is
\begin{equation}
\displaystyle\binom{G}{rK_1} = 
\displaystyle\binom{v(G)}{r} - \sum_{j\mid v(\Lambda_j) = r}N_{Ij}
\end{equation}
where indices $j$ in the summation are determined by
Lemma~\ref{lem-num-vertices}.
\end{proof}


We are interested in the question of reconstructing
a graph $G$ or some of its invariants given $\mathcal{N}(G)$.
As indicated earlier, we will assume that
the induced subgraphs $\Lambda_i; i \in [1,I]$ are not given.
We have the following relationship between Ulam's conjecture
and the $\mathcal{N}$-matrix reconstructibility.

\begin{prop}
\label{prop:equiv-un}
Ulam's conjecture is true if and only if
all graphs on three or more vertices are
$\mathcal{N}$-matrix reconstructible.
\end{prop}
\begin{proof}

\noindent 
Proof of {\em if}: by Lemma~\ref{lem-kelly1},
$\mathcal{N}(G)$ is constructed from $\mathcal{VD}(G)$.
Therefore, Ulam's conjecture is true if 
all graphs are $\mathcal{N}$-matrix reconstructible.
In fact, a graph is reconstructible if it is
$\mathcal{N}$-matrix reconstructible.\vspace{0.1in}

\noindent Proof of {\em only if}: 
this is proved by induction on the number of vertices. 
Let Ulam's conjecture be true. 
Since $N_{i1} = e(\Lambda_i)$ for all $i$, every non-empty
three vertex graph is $\mathcal{N}$-matrix reconstructible.
Now, let all graphs on at most $n$ vertices,
where $n\geq 3$, be $\mathcal{N}$-matrix reconstructible.
Let $G$ be a graph on $n+1$ vertices.
By Lemma~\ref{lem-kelly3}, the collection
$\{\mathcal{N}(G-u); u\in VG, e(G-u) > 0\}$ is unambiguously
determined by $\mathcal{N}(G)$.
The number of empty graphs in $\mathcal{VD}(G)$ is
0 or 1, and is determined by Lemma~\ref{lem-empty}.
Therefore, by induction hypothesis, $\mathcal{VD}(G)$ is uniquely
determined. Now the result follows from the assumption that
Ulam's conjecture is true.
\end{proof}

Since $\mathcal{N}(G)$ and $\mathcal{ELP}(G)$ are equivalent
by Lemma~\ref{lem-elp_eq_n}, we rephrase
Proposition~\ref{prop:equiv-un} as follows.
\begin{prop}
\label{prop:equiv-elp}
Ulam's conjecture is true if and only if all graphs
on three or more vertices are reconstructible from their
abstract edge labelled posets.
\end{prop}


We would like to point out that reconstructing
$G$ from $\mathcal{N}(G)$ or from $\mathcal{ELP}(G)$ is
not proved to be equivalent to reconstructing $G$ from
$\mathcal{VD}(G)$. This poses a difficulty. For example,
proving $\mathcal{N}$-matrix reconstructibility
of disconnected graphs is as hard as Ulam's
conjecture, although disconnected graphs are known to
be vertex reconstructible. This is proved below.

For graphs $X$ and $Y$, we use the notation $X+Y$ to
denote a graph that is a disjoint union of two graphs isomorphic to
$X$ and $Y$, respectively.
Suppose $G$ and $H$ are connected graphs having the same vertex deck. 
Consider graphs $2G = G+G$ and $2H=H+H$.

\begin{lem}
\label{lem-kelly2}
Let $F$ be a graph on fewer than $2v(G)$ vertices.
If $F$ has a component isomorphic to $G$ (in which case we write
$F=G+X$) then $\displaystyle\binom{2G}{G+X} = \displaystyle\binom{2H}{H+X}$. If $F$ has no
component isomorphic to $G$ then $\displaystyle\binom{2G}{F} = \displaystyle\binom{2H}{F}$.
\end{lem}
\begin{proof}
When $F=G+X$, $X$ must have fewer than $v(G)-1$ vertices.
Since $G$ and $H$ have identical vertex decks, by Kelly's Lemma
~\ref{lem-kelly},
$\displaystyle\binom{G}{X} = \displaystyle\binom{H}{X}$. 
Therefore, $\displaystyle\binom{2G}{G+X} = \displaystyle\binom{2H}{H+X}$.

When $F$ does not have a component isomorphic to $G$,
then if $F$ has a component on $v(G)$ vertices then
$\displaystyle\binom{2G}{F} = \displaystyle\binom{2H}{F} = 0$. 
Therefore, assume that all components of $F$ have at most
$v(G)-1$ vertices.
Any realisation of $F$ as an induced subgraph
of $2G$ is a disjoint union of graphs isomorphic to $X$
and $Y$ such that $X$ is an induced subgraph of one component
of $2G$ and $Y$ is an induced subgraph of the other component
of $2G$. Moreover, $v(X) < v(G)$ and $v(Y) < v(G)$.
Now $\displaystyle\binom{2G}{F}=\displaystyle\binom{2H}{F}$ follows from the fact that
$G$ and $H$ have identical vertex decks and Kelly's Lemma.
\end{proof}

The following corollary is an immediate consequence of the
above lemma.
\begin{cor}
\label{cor-correspondence}
Define a correspondence $f$ between $\Lambda(2G)$
and $\Lambda(2H)$ as follows.
\begin{enumerate}
\item $f(2G) = 2H$
\item $F\in \Lambda(2G)$ is not $2G$ but has a component isomorphic to $G$.
We write $F= G+X$, and set $f(F) = H+X$.
\item $F\in \Lambda(2G)$ has no component isomorphic to $G$. In this case
we set $f(F) = F$.
\end{enumerate}
The correspondence defined above is a bijection.
\end{cor}

\begin{lem}
\label{lem-2g2h}
$\mathcal{N}(2G) = \mathcal{N}(2H)$.
\end{lem}

\begin{proof}
For the bijection $f$ between non-empty induced subgraphs of $2G$ and
$2H$ that was defined in Corollary~\ref{cor-correspondence},
we show that, for any two nonisomorphic induced subgraphs
$F_1$ and $F_2$ of $2G$,

\begin{equation}
\label{eq-2g2h}
\displaystyle\binom{F_2}{F_1} = \displaystyle\binom{f(F_2)}{f(F_1)}
\end{equation}

In view of Corollary~\ref{cor-correspondence},
it is sufficient to show this when at least one of the graphs
$F_1$ and $F_2$ has a component isomorphic to $G$.
\begin{enumerate}
\item When $F_2 = 2G$, then Equation~(\ref{eq-2g2h}) follows
from Lemma~\ref{lem-kelly2} and Corollary~\ref{cor-correspondence}.

\item When $F_1 = G+X$ and $F_2 = G+Y$, and $v(X) < v(G)$ and
$v(Y) < v(G)$, we have 
$\displaystyle\binom{F_2}{F_1} = \displaystyle\binom{Y}{X} = \displaystyle\binom{H+Y}{H+X} = \displaystyle\binom{f(F_2)}{f(F_1)}$.

\item $F_2 = G+Z$, $v(Z) < v(G)$, but $F_1$ has no component
isomorphic to $G$. In this case, any realisation of $F_1$ as
an induced subgraph of $F_2$ may be represented (possibly in many ways)
as $F_1 = X+Y$ where $X$ is an induced proper subgraph of the component
$G$ of $F_2$ and $Y$ is an induced subgraph of $Z$.
Moreover, $v(X) < v(G)$ and
$v(Y) < v(G)$. Since $\displaystyle\binom{G}{X} = \displaystyle\binom{H}{X}$, we have
$\displaystyle\binom{G+Z}{F_1} = \displaystyle\binom{H+Z}{f(F_1)}$.
Note that the actual value of
$\displaystyle\binom{F_2}{F_1}$ may be written by considering all possible
ways of realising $F_1$ as an induced subgraph of $G+Z$.
\end{enumerate}

Thus we have shown Equation~(\ref{eq-2g2h}) for arbitrary
non-empty induced subgraphs of $2G$, which implies the result.
\end{proof}

\begin{thm}
\label{thm-equiv-disc}
Ulam's conjecture is true if and only if
disconnected graphs on three or more vertices are
$\mathcal{N}$-matrix reconstructible.
\end{thm}
\begin{proof} The {\em only if} part follows from
Proposition~\ref{prop:equiv-un}.
The {\em if} part is proved by contra\-diction.
Suppose $G$ and $H$ are connected nonisomorphic graphs with
the same vertex deck, that is, they are a counter example to
Ulam's conjecture. Then $2G$ and $2H$ are nonisomorphic
but $\mathcal{N}(2G) = \mathcal{N}(2H)$ by Lemma~\ref{lem-2g2h}.
Therefore, $2G$ and $2H$ are disconnected graphs that are
not $\mathcal{N}$-matrix reconstructible.
\end{proof}


The following result is proved along the lines of Lemma~\ref{lem-2g2h}.

\begin{thm}
Ulam's conjecture is true if and only if the edge labelled poset
of each graph has only the trivial automorphism.
\end{thm}

\begin{proof}
The proof of {\em only if} is done by contradiction.
Suppose that $\mathcal{ELP}(G)$ has a nontrivial automorphism $\sigma$.
Then there are nonisomorphic induced subgraphs $\Lambda_i$ and $\Lambda_j$ 
of $G$ such that $\sigma(\Lambda_i) = \Lambda_j$. The downsets
(or the edge labelled posets) of $\Lambda_i$ and $\Lambda_j$ must be
isomorphic. Therefore, by Proposition~\ref{prop:equiv-elp},
there is a counter example to Ulam's conjecture.

The proof of {\em if} is also done by contradiction.
Suppose that Ulam's conjecture is false, and $G$ and $H$
are connected nonisomorphic graphs having identical vertex decks.
We show that $\mathcal{ELP}(G+H)$ has a nontrivial
automorphism. 
Define a bijective map 
$\sigma : \Lambda(G+H)\rightarrow  \Lambda(G+H)$ as follows. 
\begin{enumerate}
\item The graph $G+H$ is mapped to itself.
\item If $\Lambda_i\in \Lambda(G+H)$ has a component isomorphic to $G$,
then denote $\Lambda_i$ by $G+X$, where $X$ is a proper subgraph
of the component isomorphic to $H$. In this case,
set $\sigma(G+X) = H+X$.
\item If $\Lambda_i$ is $H+X$, where $X$ is a proper subgraph
of the component isomorphic to $G$, then set
$\sigma(H+X) = G+X$.
\item For all other graphs $\Lambda_i \in \Lambda(G+H)$,
$\sigma(\Lambda_i) = \Lambda_i$.
\end{enumerate}

We now show that $\sigma$ is an automorphism of
$\mathcal{ELP}(G+H)$. That is, we show that
$\displaystyle\binom{\Lambda_i}{\Lambda_j} = 
\displaystyle\binom{\sigma(\Lambda_i)}{\sigma(\Lambda_j)}$ for any
two graphs $\Lambda_i $ and $\Lambda_j$ in $\Lambda(G+H)$.

We have to consider only the case in which at least one of
$\Lambda_i $ and $\Lambda_j$ has a component isomorphic to $G$ or
$H$, and $v(\Lambda_j) \leq v(\Lambda_i)$.

\begin{enumerate}

\item $\Lambda_j = G+X$ and $\Lambda_i = G+H$.
In this case, \\
$\displaystyle\binom{G+H}{G+X} = \displaystyle\binom{H}{X} = \displaystyle\binom{G}{X} = 
\displaystyle\binom{H+G}{H+X} = \displaystyle\binom{\sigma(G+H)}{\sigma(G+X)}$.

\item $\Lambda_j = G+X$ and $\Lambda_i = G+Y$ and $v(Y) < v(G)=v(H)$.
In this case, \\ 
$\displaystyle\binom{G+Y}{G+X} = \displaystyle\binom{Y}{X} = \displaystyle\binom{H+Y}{H+X} = 
\displaystyle\binom{\sigma(G+Y)}{\sigma(G+X)}$.

\item $\Lambda_j = G+X$ and $\Lambda_j = H+Y$ and $Y\ncong G$.
In this case, \\
$\displaystyle\binom{H+Y}{G+X} = \displaystyle\binom{G+Y}{H+X} = 0$.
 
\item $\Lambda_j = G+X$ and $\Lambda_i $ has no component isomorphic to
$G$ or $H$. In this case,\\
$\displaystyle\binom{\Lambda_i}{G+X} = \displaystyle\binom{\Lambda_i}{H+X} = 0$.

\item $\Lambda_j$ has no component isomorphic to $G$ or $H$,
and $\Lambda_i = G+H$. This is trivial since 
$\sigma(\Lambda_j) = \Lambda_j$ and $\sigma(G+H) = G+H$

\item $\Lambda_j$ has no component isomorphic to $G$ or $H$
and $\Lambda_i = G+X$, where $v(X) < v(G)=v(H)$. In this case,
a realisation of $\Lambda_j$ as an induced subgraph of $G+X$
may be written as $\Lambda_j = Y+Z$, where $Y$ is an induced subgraph
of $G$ and $Z$ is an induced subgraph of $X$. Since,
$\displaystyle\binom{G}{Y}=\displaystyle\binom{H}{Y}$, the number of such realisations is
$\displaystyle\binom{G}{Y}\displaystyle\binom{X}{Z} = \displaystyle\binom{H}{Y}\displaystyle\binom{X}{Z}$. By summing
over all possible ways of realising $\Lambda_j$ as an induced
subgraph of $G+X$, we get $\displaystyle\binom{G+X}{\Lambda_j} = \displaystyle\binom{H+X}{\Lambda_j}$.

\item All the above arguments are valid when $G$ and $H$ are
interchanged.
\end{enumerate}

Thus we have constructed a non-trivial automorphism, completing the
{\em if} part.
\end{proof}


We conclude this section with a result on trees.
\begin{thm}
\label{thm-trees}
Trees and forests are $\mathcal{N}$-matrix reconstructible.
\end{thm}
\begin{proof}
The class of simple acyclic graphs is closed under
vertex deletion. Therefore, we can use the method
in the proof of Proposition~\ref{prop:equiv-un}.
Let $T$ be a tree or a non-empty forest on three or more vertices.
We prove by induction on $v(T)$ that $T$ is uniquely
reconstructible from $\mathcal{N}(T)$. The base case is
$v(T) = 3$. All graphs on 3 vertices are
$\mathcal{N}$-matrix reconstructible by Lemma~\ref{lem-num-vertices}.
Suppose each acyclic graphs on at most $k$
can be recognised and reconstructed from its
$\mathcal{N}$-matrix. Let $v(T) = k+1$.
By Lemma~\ref{lem-kelly3}, the collection
$\{\mathcal{N}(T-u)| u \in VT\}$ is unambiguously
determined. Then by induction hypothesis, $T-u$ are
determined (along with their multiplicities). The
subgraphs in the vertex deck that are not determined by
Lemma~\ref{lem-kelly3} are the ones having no edges.
Since Ulam's conjecture has been proved for
trees and disconnected simple graphs in \cite{kelly1957},
$T$ is $\mathcal{N}$-matrix reconstructible.
\end{proof}

\noindent {\bf Remark} If Ulam's conjecture is true
for a class of graphs that is closed under vertex
deletion, then the class is also $\mathcal{N}$-matrix
reconstructible.


\section{Tutte-Kocay theory on the $\mathcal{N}$-matrix.}
\label{sec-nmatrix}

In this section we will compute several invariants of a graph
$G$ from its $\mathcal{N}$-matrix.
The invariants include the number of spanning trees,
the number of spanning unicyclic subgraphs containing 
a cycle of specified length, the characteristic polynomial
and the rank polynomial.

\noindent {\bf An outline of the proof.}
First we outline how the above mentioned invariants are
calculated from the vertex deck using Kocay's Lemma.
\begin{enumerate}
\item Suppose the graphs $F_1, F_2, \ldots , F_k$ satisfy
$\sum_i v(F_i) = v(G)$ and $v(F_i)< v(G) \forall i$.
Kocay's Lemma then gives the number
of disconnected spanning subgraphs having components isomorphic to
$F_1, F_2, \ldots , F_k$.
\item Kacay's lemma is then applied to 
$F_1 = F_2 = \ldots = F_k = K_2$, where $k=v(G)-1$.
Since disconnected spanning subgraphs of each type are counted
in the first step, we can now count the number of spanning trees.
\item The second step is repeated with $k=v(G)$. Since
the number of spanning trees and disconnected spanning subgraphs
of each type are known from the first two steps, we can now
count the number of hamiltonian cycles.
\item Once the above three steps are completed, many other invariants,
such as the characteristic polynomial, rank polynomial, etc. are
easily computed.
\end{enumerate}

The procedure outlined above cannot be implemented 
on the $\mathcal{N}$-matrix in a straight forward manner.
We do not know all the induced proper
subgraphs. But we observe that the above procedure essentially
reduces counting certain spanning subgraphs to counting
them on vertex proper subgraphs. It turns out that we
do not really need the number of vertex proper subgraphs
of each type.
We only need to know the `cycle structure', that is,
$\psi_i(\Lambda_j)$ for each $i \leq v(\Lambda_j)$, for each $j < I$.
Next we outline the strategy to construct the cycle structure.

Suppose $X, Y, \ldots$ is a list of some graph invariants that are
either polynomials or numbers, for example, the number of hamiltonian
cycles in a graph or the chromatic polynomial of a graph.
We say that an invariant $Z$ can be reduced to invariants
$X, Y, \ldots $ (or $Z$ has a reduction on the $\mathcal{N}$-matrix)
if for each graph $G$ having a non-empty edge set,\
\begin{enumerate}
\item $Z(G)$ can be written as $Z(G) = \Theta(X(G_U), Y(G_V), ...)$
where $\Theta(x,y,...)$ is a polynomial in $x,y, \ldots$,
and $U, V, \dots $ are proper subsets of $VG$.
\item the coefficient of each term in the polynomial can be
computed from $\mathcal{N}(G)$.
\end{enumerate}
Proving an identity that gives a reduction of an invariant
$Z$ as in the above equation is not in itself
sufficient to claim that $Z$ is $\mathcal{N}$-matrix
reconstructible. If the invariants $X_1, X_2, \ldots, X_k$ appear on the
RHS of the above equation, then it is essential to
show that the invariants $X_1, X_2, \ldots, X_k$ themselves
can be reduced to $X_1, X_2, \ldots, X_k$.
The reconstructibility of $Z$ and $X_1, X_2, \ldots, X_k$
from the $\mathcal{N}$-matrix is then proved by 
induction on $v(G)$. That is possible because of the
requirement that the sets $U,V,\ldots $ are proper subsets
of $VG$. It is worth noting that the chromatic polynomial
computation given in Biggs \cite{biggs1993} essentially
follows a similar style.
In several lemmas that precede the main theorem, we will
prove identities of the form $Z(G) = \Theta(X(G_U), Y(G_V), ...)$.
It will become clear that in the end all invariants computed
here will reduce to the cycle structure of proper subgraphs.

\begin{lem}
\label{lem-kelly-cycles}
For $i < v(G)$, the number of cycles of length $i$ in $G$ has a
reduction on $\mathcal{N}(G)$ given by
\begin{equation}
\psi_i(G) = \frac{1}{v(G)-i}\sum_{u\in VG}\psi_i(G-u)
= \frac{1}{v(G)-i}\sum_{j\mid v(\Lambda_j)=v(G)-1}N_{Ij}\psi_i(\Lambda_j)
\end{equation}
\end{lem}
\begin{proof} This immediately follows from Kelly's Lemma ~\ref{lem-kelly}
and Lemma~\ref{lem-kelly3}.
\end{proof}

\begin{defn}
\label{def-cycle-cover}
Let $X$ be a subset of $VG$. Let $A\equiv (a_i)_{i=1}^k$
be a sequence in $[2,|X|]$.
A $k$-tuple of cycles in $G_X$,
corresponding to the sequence $A$, is a $k$-tuple
of cycles in $G_X$ such that the $k$ cycles 
have lengths $a_1, a_2, \ldots, a_k$, respectively.
Additionally, if the cycles in the $k$-tuple span the
set $X$, then it is called a {\em spanning cycle cover}
of $G_X$, corresponding to the sequence $A$.
The number of $k$-tuples of cycles in $G_X$, corresponding
to the sequence $A$, is denoted by $p(A\rightarrow G_X)$.
The number of spanning cycle covers of $G_X$, corresponding
to the sequence $A$, is denoted by $c(A\rightarrow G_X)$.
\end{defn}

\begin{lem}
\label{lem-p2c}
\begin{equation}
\label{eq-p2c}
p(A \rightarrow G_X) = \sum_{Y\subseteq X}c(A\rightarrow G_Y)
\end{equation}
\end{lem}

\begin{proof}
\begin{equation}
\label{eq-p2c-proof}
\begin{split}
p(A \rightarrow G_X) & = \prod_{j=1}^{k}\psi_{a_j}(G_X)\\  
 & = |\{(F_1,F_2,\ldots,F_{k})\mid (\forall j\in[1,k])(F_j\subseteq G_X,\,
          F_j\cong C_{a_j})\}| \\
 & = \sum_{Y\subseteq X}c(A\rightarrow G_Y)
\end{split}
\end{equation}
Thus we have grouped together the $k$-tuples of cycles
in groups that span each subset of $X$. This is
essentially the idea of Kocay's Lemma~\ref{lem-kocay}.
\end{proof}

\begin{lem}
\label{lem-exp-p}
If $A \equiv (a_i)_{i = 1}^k $ is a sequence in $[2,v(G)-1]$ then
$p(A\rightarrow G)$ has a reduction on $\mathcal{N}(G)$ given by
\begin{equation}
p(A \rightarrow G) =\prod_{i=1}^{k}\psi_{a_i}(G)
= \prod_{i=1}^{k}
\left(\frac{1}{v(G)-a_i}
\sum_{j\mid v(\Lambda_j)=v(G)-1}N_{Ij}\psi_{a_i}(\Lambda_j)\right)
\end{equation}
\end{lem}
\begin{proof}
This follows from the definition of $p(A\rightarrow G)$ and
Lemma~\ref{lem-kelly-cycles}.
\end{proof}

\begin{lem}
\label{lem-exp-c}
If $A \equiv (a_i)_{i = 1}^k $ is a sequence in $[2,v(G)-1]$ then
$c(A\rightarrow G)$ has a reduction on $\mathcal{N}(G)$.
\end{lem}

\begin{proof}
By M\"{o}bius inversion of Equation~(\ref{eq-p2c}), we write
\begin{equation} 
c(A\rightarrow G_X) = \sum_{Y\subseteq X}(-1)^{|X\backslash Y|}
p(A\rightarrow G_Y)
\end{equation}
which implies
\begin{equation}
c(A\rightarrow G) = 
\sum_{j=1}^I(-1)^{v(G)- v(\Lambda_j)}N_{Ij}p(A \rightarrow \Lambda_j)
\end{equation}
By Lemma~\ref{lem-exp-p}, the RHS of this equation has a
reduction on $\mathcal{N}(G)$. Therefore, $c(A\rightarrow G)$ has
a reduction on $\mathcal{N}(G)$.
\end{proof}

The following definition restricts the spanning cycle covers
of Definition~\ref{def-cycle-cover} to connected
spanning cycle covers.
\begin{defn}
\label{def-connected-spanning-cycle-cover}
Let $X$ be a subset of $VG$.
Let $A\equiv (a_i)_{i=1}^k$ be a sequence in $[2,|X|]$.
A {\em connected spanning cycle cover} of $G_X$,
corresponding to the sequence $A$, is a $k$ tuple
of cycles in $G_X$ such that the $k$ cycles 
have lengths $a_1, a_2, \ldots, a_k$, respectively,
and together they constitute a connected subgraph spanning 
the set $X$. More formally, it is a $k$-tuple
$(F_1,F_2,\ldots,F_{k})$ such that
$(\forall j\in [1,k])(F_j\subseteq G_X, F_j\cong C_{a_j})$,
$\cup_{j=1}^k VF_j = X$, and $\cup_{j=1}^kF_j$ is connected.
The number of connected spanning cycle covers of $G_X$,
corresponding to the sequence $A$, is denoted by
$con(A\rightarrow G_X)$. The {\em disconnected
spanning cycle covers} are defined similarly, and their
number, corresponding to a sequence $A$, is denoted by 
$discon(A\rightarrow G_X)$.
\end{defn}

Let $A\equiv (A_i)_{i=1}^l$ be a list of $l$ non-increasing
sequences such that 
$A_i \equiv (a_{ij})_{j=1}^{k_i}$; $a_{ij}\in [2,v(G)]$.
Let $B \equiv (b_i)_{i=1}^l$ be a sequence in $[2,v(G)]$.
Let $m \leq n$. We now define quantities $Q_m(A,B\rightarrow G)$
and $T_p(A,B\rightarrow G)$ that are based on connected spanning
cycle covers as follows.

\begin{equation}
\begin{split}
Q_m(A,B\rightarrow G)
& =  \sum_{\substack {S\subseteq V(G)\\
                       |S|=m}}
        \prod_{i = 1}^{l}
        \left(
        \sum_{\substack{X\subseteq S \\
                        |X|=b_i}}
        con(A_i\rightarrow G_{X}) 
        \right)\\
& =  \sum_{\substack{S\subseteq V(G)\\|S|=m}}
        \left(
        \sum_{\substack{(X_1,X_2,\ldots,X_l)\mid \\
                        \bigcup_{j=1}^{l}X_j\subseteq S\\
                        |X_j|=b_j\forall j}}
        \left(
        \prod_{i = 1}^{l}con(A_i\rightarrow G_{X_i}) 
        \right)
        \right)\\
& = \sum_{p\leq m} T_p(A,B\rightarrow G)\displaystyle\binom{v(G)-p}{m-p}
\end{split}
\end{equation}
where
\begin{equation}
\begin{split}
 T_p(A,B\rightarrow G) & =
        \sum_{\substack{(X_1,X_2,\ldots,X_l)\mid \\
                        \cup_{j=1}^{l}X_j\subseteq V(G)\\
                        |\cup_{j=1}^{l}X_j| = p \\
                        |X_j|=b_j\forall j}
        }
        \left(
        \prod_{i = 1}^{l}con(A_i\rightarrow G_{X_i}) 
        \right)
\end{split}
\end{equation}

Solving the system of equations for
$T_m(A,B\rightarrow G)$, we can write
\begin{equation}
\label{eq-q2tp}
\begin{split}
T_m(A,B\rightarrow G) & =
\sum_{p\leq m} (-1)^{m-p}
\displaystyle\binom{v(G)-p}{m-p}Q_p(A,B\rightarrow G)
\end{split}
\end{equation}
When $m=v(G)$, this is simply
\begin{equation}
\begin{split}
T_{v(G)}(A,B\rightarrow G) & =
\sum_{p\leq v(G)} (-1)^{v(G)-p} Q_p(A,B\rightarrow G)
\end{split}
\end{equation}

Note that if $m < max_{i,j}(a_{ij})$ for some $i,j$ then
$T_m(A,B\rightarrow G)$ and $Q_m(A,B\rightarrow G)$
are both 0.

\begin{lem}
\label{lem-exp-q-tp}
If $A_i; i \in [1,l]$ are sequences in $[2,v(G)-1]$,
and $B\equiv (b_i)_{i=1}^l$ are sequences in $[2,v(G)-1]$
then $Q_m(A,B\rightarrow G)$ and $T_m(A,B\rightarrow G)$
have reductions on the $\mathcal{N}$-matrix
for each $m \leq v(G)$.
\end{lem}

\begin{proof}
We write $Q_m(A,B\rightarrow G)$ in terms of $\Lambda_j$
and the entries of $\mathcal{N}(G)$.
\begin{equation}
\label{eq-exp-q}
\begin{split}
Q_m(A,B\rightarrow G)
& =  \sum_{\substack {S\subseteq V(G)\\
                       |S|=m}}
        \prod_{i = 1}^{l}
        \left(
        \sum_{\substack{X\subseteq S \\
                        |X|=b_i}}
        con(A_i\rightarrow G_{X}) 
        \right)\\
&= \sum_{p\mid v(\Lambda_p)=m}N_{Ip}
        \prod_{i = 1}^{l}
        \left(
        \sum_{j\mid v(\Lambda_j)=b_i}N_{pj}con(A_i\rightarrow \Lambda_j)
        \right)
\end{split}
\end{equation}

Since $b_i<v(G)$, Equations~(\ref{eq-exp-q})
reduce $Q_m(A,B\rightarrow G)$ to the invariants
$con(A_i\rightarrow \Lambda_j)$ for $m\leq v(G)$.
Therefore, by Equation~(\ref{eq-q2tp}), $T_m(A,B\rightarrow G)$
are also reduced to the invariants $con(A_i\rightarrow \Lambda_j)$
for each $m\leq v(G)$.
Note that if $a_{ik} > v(\Lambda_j)$ for some $k$ then 
$con(A_i\rightarrow \Lambda_j)=0$.
\end{proof}

\begin{lem}
\label{lem-con}
If $A\equiv (a_i)_{i=1}^k$ is a sequence in $[2,v(G)-1]$
then $con(A\rightarrow G)$ has a reduction on $\mathcal{N}(G)$.
\end{lem}

\begin{proof} The idea of the proof is similar to Kocay's Lemma.
First $c(A\rightarrow G)$ is written as 
$con(A\rightarrow G) + discon(A\rightarrow G)$.
Then $discon(A\rightarrow G)$ is expressed in terms
of $con(A_i \rightarrow G_{X_i})$ where 
$X_i$ are proper subsets of $VG$, and $A_i$ are certain subsequences
of $A$.
Then $discon(A\rightarrow G)$ are related to 
$T_{v(G)}(A^P, B\rightarrow G)$ for certain subsequences $A^P$ of $A$ 
and certain sequences $B$ constructed from appropriate partitions of $VG$.
Since the reductions of $c(A\rightarrow G)$ and
$T_{v(G)}(A^P, B\rightarrow G)$ have already been obtained,
we get a reduction of $con(A\rightarrow G)$. 

Let $\mathcal{P}_q(N_k)$ be the set of all partitions of
$N_k$ in $q$ parts. 
A partition $P$ in $\mathcal{P}_q(N_k)$,
is denoted by $P = \{N_k^1,N_k^2,\ldots ,N_k^q\}$.
Consider an arbitrary $k$ tuple $(F_1,F_2,\ldots,F_{k})$ 
such that
$(\forall i\in [1,k])(F_i\subseteq G,F_i\cong C_{a_i})$.
It defines a partition $P\in \mathcal{P}_q(N_k)$
so that $h$ and $j$ are in the same part of $P$ if
and only if $F_h$ and $F_j$ are subgraphs of the same
connected component of $\cup_{i=1}^kF_i$. We denote the contribution
to $c(A\rightarrow G)$ from such tuples by 
$c_P(A\rightarrow G)$, and write
\begin{equation}
c(A\rightarrow G) = \sum_{q}\sum_{P\in \mathcal{P}_q}c_P(A\rightarrow G)
\end{equation}
Let the set of solutions to the equation $\sum_{i=1}^qb_i = v(G)$
be $\mathcal{B}(q)$. We can then write
\begin{equation}
c(A\rightarrow G) =
\sum_{q}
\sum_{P\in \mathcal{P}_q}
\sum_{B\in \mathcal{B}(q)}
c_{PB}(A\rightarrow G)
\end{equation}
where one more suffix $B$ in the summand is used to denote those
tuples for which the connected component corresponding
to part $N_k^i$ has $b_i$ vertices, for $i\in [1,q]$.
Now expanding the summand in terms of $con(\ldots)$, we get
\begin{equation}
\begin{split}
c(A\rightarrow G) =
\sum_{q}
\sum_{P\in \mathcal{P}_q}
\sum_{B\in \mathcal{B}(q)}
\sum_{\substack{(X_1,X_2,\ldots,X_q)\mid \\
                \cup_{j=1}^{q}X_j= V(G)\\
                |X_j|=b_j\forall j}
     }
\prod_{i=1}^q con(A_i\rightarrow G_{X_i})
\end{split}
\end{equation}
where $A_i$ is the subsequence of $A$ with indexing set $N_k^i$.
Innermost summation and product are now replaced by $T_{v(G)}(A^P,B\rightarrow G)$,
so
\begin{equation}
\label{eq-pqb}
\begin{split}
c(A\rightarrow G) =
\sum_{q}
\sum_{P\in \mathcal{P}_q}
\sum_{B\in \mathcal{B}(q)}
T_{v(G)}(A^P,B\rightarrow G)
\end{split}
\end{equation}
where $A^P$ is the collection of subsequences
$A_i$ of $A$; $i\in [1,q]$, corresponding to the partition $P$.

By Lemma~\ref{lem-exp-c}, the LHS of Equation~(\ref{eq-pqb})
has a reduction on $\mathcal{N}(G)$.
By Lemma~\ref{lem-exp-q-tp}, each term on the RHS,
except the term $con(A\rightarrow G)$, which corresponds to $q = 1$,
has a reduction on $\mathcal{N}(G)$.
This proves that $con(A\rightarrow G)$ has a reduction on $\mathcal{N}(G)$.
\end{proof}

\begin{cor}
\label{cor-trees}
Let $tr(G)$ be the number of spanning trees in $G$, and let
$uni(G,r)$ be the number of spanning unicyclic subgraphs of $G$
containing an $r$-cycle, respectively. Then, $tr(G)$ and 
$uni(G,r); r\in [3,v(G)]$ have reductions on $\mathcal{N}(G)$.

\end{cor}

\begin{proof}
Define $A\equiv (a_i)_{i=1}^{v(G)-1}$ such that
$a_i=2 \,\forall \, i\in [1,v(G)-1]$.
By Lemma~\ref{lem-con}, $con(A\rightarrow G)$
has a reduction on $\mathcal{N}(G)$. But
$con(A\rightarrow G) = (v(G)-1)!tr(G)$. Therefore,
$tr(G)$ has a reduction $\mathcal{N}(G)$.

To reduce $uni(G,r);r\in [3,v(G)-1]$,
define $A\equiv(a_j)_{j=1}^{v(G)-r+1}$ such that 
$a_1 = r$, and $a_j = 2 \,\forall \, j\in [2,v(G)-r+1]$.
Again, $con(A\rightarrow G) = (v(G)-r)!uni(G,r)$
has a reduction $\mathcal{N}(G)$ by Lemma~\ref{lem-con}.
Thus $uni(G,r)$ also  has a reduction on $\mathcal{N}(G)$.

To reduce the number of hamiltonian cycles,
let $A\equiv(a_i)_{i=1}^{v(G)}; a_i = 2\,\forall \,i\in [1,v(G)]$.
We have,
$con(A\rightarrow G) = (v(G)-1)!S(v(G),v(G)-1)tr(G) + 
\sum_{i=3}^{v(G)-1}v(G)!uni(G,i) + v(G)!ham(G)$,
where $S(v(G),v(G)-1)$ is the Sterling number of the second kind
computed for $(v(G),v(G)-1)$.
Therefore, $ham(G)$, has a reduction $\mathcal{N}(G)$.
\end{proof}

\begin{lem}
\label{lem-ci}
The coefficients $c_i(G)$ of the characteristic polynomial
of $G$ are given by $c_0(G) = 1$ and 
\begin{equation}
c_i(G) = \frac{1}{v(G)-i}\sum_{j\mid v(\Lambda_j)=v(G)-1}N_{Ij}c_i(\Lambda_j) 
\,\,\text{for}\,\, 0 < i < v(G).
\end{equation}
\end{lem}

\begin{proof} By Lemma~\ref{lem-derivative} we write
$P'(G;\lambda) = \sum_{u\in VG}P(G-u;\lambda)$.
Equating identical powers of $\lambda $ on the two sides,
we get the result.
\end{proof}

\begin{cor}
\label{cor-elementary}
If $F$ is an elementary graph on $v(G)$ vertices
then $\ncra{G}{F}$ has a reduction on $\mathcal{N}(G)$.
\end{cor}

\begin{proof} 
The case when $F$ is a hamiltonian cycle is handled in
Corollary ~\ref{cor-trees}. 
If $F$ is not a cycle,
define a sequence $A\equiv (a_i)_{i=1}^k,\,
a_i \in [2,v(G)-1]$ such that $\sum_{i=1}^k a_i = v(G)$.
It is uniquely associated with an elementary graph
$F\in L_{v(G)}$, so that the components of $F$ are cycles of
length $a_i$, or $K_2$ if $a_i=2$. 
Now $c(A\rightarrow G) = c(A\rightarrow F)\ncra{G}{F}$,
and $c(A\rightarrow F)$ depends only on the multiplicity
of each cycle length in $F$. By Lemma~\ref{lem-exp-c},
$c(A\rightarrow G)$ has a reduction on the
$\mathcal{N}$-matrix. Therefore, $\ncra{G}{F}$ has a
reduction on the $\mathcal{N}$-matrix.
\end{proof}

The following chart shows how various invariants were
reduced to other invariants. For example, it shows
that $con(A\rightarrow G); \, 2 \leq a_i \leq v(G)$ can be
reduced to computing 
$ham(G)$, $con(A\rightarrow G); \, 2 \leq a_i < v(G)$
and $p(A\rightarrow G); \, 2 \leq a_i < v(G)$. It is clear
from the diagram that computing 
$con(A\rightarrow G); \, 2 \leq a_i \leq v(G)$ and $ham(G)$
reduces to computing the same invariants for induced proper
subgraphs.

\newpage
\begin{figure}[ht]
\label{fig-dependency}
\input{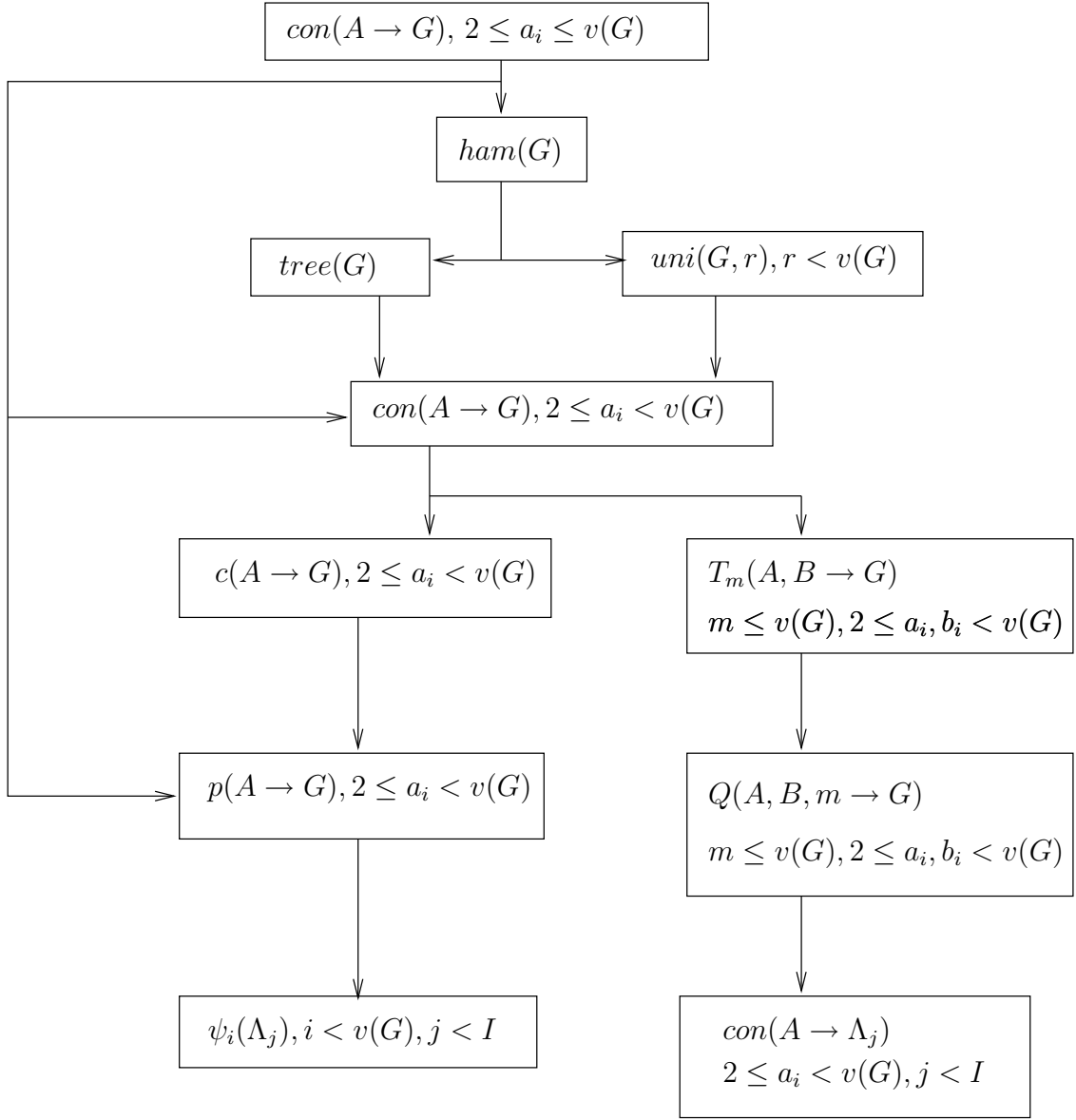}
\vspace{0.2in}
\caption[]{A summary of invariant reductions}
\end{figure}

This makes the reconstructibility of several invariants obvious.
\begin{thm}
\label{thm-everything}
Suppose that $G$ is a simple finite graph, and we are
given $\mathcal{N}(G)$.
\begin{enumerate}
\item Let $A\equiv (a_i)_{i=1}^k$ be a sequence
such that $a_i \in [2,v(G)]$. Then 
$con(A\rightarrow G)$ and $ham(G)$ are reconstructible from
$\mathcal{N}(G)$.
\item 
If $A_i; i \in [1,l]$ are sequences in $[2,v(G)]$,
and $B\equiv (b_i)_{i=1}^l$ is a sequence in $[2,v(G)]$
then $Q_m(A,B\rightarrow G)$ and $T_m(A,B\rightarrow G)$
are reconstructible from $\mathcal{N}(G)$ for each $m \leq v(G)$.
\item the number of spanning trees in $G$,
the number of cycles of length $i$, for $3 \leq i \leq v(G)$,
the number of unicyclic subgraphs containing a cycle of
length $i$, for each $i \in [3,v(G)]$, and the
characteristic polynomial $P(G;\lambda)$ are all reconstructible
from $\mathcal{N}(G)$.
\end{enumerate}
\end{thm}

\begin{proof}
We prove the first item by induction on $v(G)$.
The base case is $v(G) = 2$. In this case
$G = K_2$. Now suppose that  $con(A\rightarrow G)$ and $ham(G)$
are reconstructible from $\mathcal{N}(G)$ when $v(G) < s$
for an arbitrary sequence $A\equiv (a_i)_{i=1}^k$ of integers
in  $[2,v(G)]$. Now let $v(G) = s$. Lemmas and Corollaries
~\ref{lem-kelly-cycles} to ~\ref{cor-trees} imply that
computations of $con(A\rightarrow G)$ and $ham(G)$
reduce to computations on induced proper subgraphs of $G$,
thus completing the induction step, and the proof of the
first item.

Since all other intermediate invariants $T_m(\ldots)$,
$Q_m(\ldots)$, the number of spanning trees, the number of
unicyclic graphs having a cycle of a specified length,
number of cycles of each length, number of elementary spanning
graphs of each type, etc. have been reduced to computations
of invariants $con(A\rightarrow \Lambda_j); j < I$ and
$\psi_i(\Lambda_j); i < v(G), j < I$, the remaining parts of the
theorem follow.

\end{proof}

There is another way of proving the $\mathcal{N}$-matrix reconstructibility
of the characteristic poly\-nomial. From $\mathcal{N}(G)$,
it is possible to construct $\mathcal{N}(\bar{G})$, and then
invoke the result of Hagos \cite{hagos2000} in the induction
step. Hagos proved that the pair $(P(G;\lambda),P(\bar{G};\lambda))$
can be reconstructed from the collection
$\{(P(G-u;\lambda),P(\bar{G}-u;\lambda)); u\in VG\}$. We skip
the details of this argument.
The proof presented here counts many other invariants.
It is likely that the deck of pairs of polynomials
considered by Hagos contains enough information for
counting hamiltonian cycles and spanning trees.

Now we count the subgraphs with a given number
of components, and a given number of edges in each component,
and use it to compute the rank polynomial.

Let $\mathcal{G}(p,l,(p_i,q_i)_{i=1}^l)$ be the family of
graphs with $p$ vertices and $l$ components, such that
the $i$'th component has $p_i$ vertices and $q_i$
edges for $i \in [1,l]$. So, $\sum_i p_i = p$. We
also assume that $p_i \geq p_j$ whenever $i < j$. 
By extending the notation $\ncra{G}{F}$ defined earlier,
we denote by $\ncra{G}{\mathcal{G}(p,l,(p_i,q_i)_{i=1}^l)}$
the number of subgraphs of $G$ that belong to the
family $\mathcal{G}(p,l,(p_i,q_i)_{i=1}^l)$.

\begin{lem}
\label{lem-kedge}
The number of connected spanning subgraphs of $G$ with
$k$ edges, that is, \\
$\ncra{G}{\mathcal{G}(v(G),1,(v(G),k))}$,
is reconstructible from $\mathcal{N}(G)$ for all $k$.
\end{lem}
\begin{proof} When $k < v(G)-1$,
$\ncra{G}{\mathcal{G}(v(G),1,(v(G),k))}$ is 0.
When $k\in [v(G)-1, e(G)]$, we prove the result
by induction on $k$. The base case $k = v(G)-1$,
which corresponds to the number of spanning trees,
was proved in Theorem~\ref{thm-everything}.
Let the claim be true for
all $k\in [v(G)-1,q-1]$. To prove the claim for
$k=q$, define $A\equiv (a_i)_{i=1}^q$ such that
$a_i=2$ for all $i\leq q$. We can write
\begin{equation}
\label{eq-kedge-con}
con(A\rightarrow G) = \sum_{i=n-1}^q i!S(q,i)
\ncra{G}{\mathcal{G}(v(G),1,(v(G),i))}
\end{equation}
In Theorem~\ref{thm-everything}, $con(A\rightarrow G)$
was shown to be reconstructible from $\mathcal{N}(G)$.
By the induction hypothesis, all terms on the RHS, except
$\ncra{G}{\mathcal{G}(v(G),1,(v(G),q))}$,
are known. Solving Equation~(\ref{eq-kedge-con}) for
$\ncra{G}{\mathcal{G}(v(G),1,(v(G),q))}$ completes the induction
step and the proof.
\end{proof}

\begin{lem}
\label{lem-lcompo}
$\ncra{G}{\mathcal{G}(v(G),l,(n_i,m_i)_{i=1}^l)}$,
where $\sum_{i=1}^ln_i = v(G)$ and $m_i > 0$ for all $i$,
is reconstructible from $\mathcal{N}(G)$.
\end{lem}
\begin{proof} 
Let $A\equiv (A_i)_{i=1}^l$, where $A_i\equiv (a_{ij})_{j=1}^{m_i}$;
$a_{ij} = 2\,\forall i,j$, and $B\equiv (n_i)_{i=1}^l$.
By Theorem~\ref{thm-everything}, 
$T_{v(G)}(A,B\rightarrow G)$ is $\mathcal{N}$-matrix 
reconstructible. 
We first express $T_{v(G)}(A,B\rightarrow G)$
in terms of the subgraphs to be counted, and then
count the subgraphs by induction.

\begin{equation}
\begin{split}
 T_{v(G)}(A,B\rightarrow G) & =
        \sum_{\substack{(X_1,X_2,\ldots,X_l)\mid \\
                        \cup_{j=1}^{l}X_j= V(G)\\
                        |X_j|=n_j\forall j}
        }
        \prod_{i = 1}^{l}con(A_i\rightarrow G_{X_i}) \\
 & =    \sum_{\substack{(X_1,X_2,\ldots,X_l)\mid \\
                        \cup_{j=1}^{l}X_j= V(G)\\
                        |X_j|=n_j\forall j}
        }
        \prod_{i = 1}^{l}\left(
\sum_{q_i = n_i-1}^{m_i} q_i!S(m_i,q_i)
\ncra{G_{X_i}}{\mathcal{G}(n_i,1,(n_i,q_i))}
\right)\\
 & =    \sum_{\substack{(X_1,X_2,\ldots,X_l)\mid \\
                        \cup_{j=1}^{l}X_j= V(G)\\
                        |X_j|=n_j\forall j}
        }
        \sum_{\substack{(q_1,q_2,\ldots,q_l)\mid \\
                        n_j-1\leq q_j \leq m_j \forall j}
        }
        \prod_{i = 1}^{l}\left(
q_i!S(m_i,q_i)\ncra{G_{X_i}}{\mathcal{G}(n_i,1,(n_i,q_i))}
\right)\\
 & =    \sum_{\substack{(q_1,q_2,\ldots,q_l)\mid \\
                        n_j-1\leq q_j \leq m_j \forall j}
        }
        \sum_{\substack{(X_1,X_2,\ldots,X_l)\mid \\
                        \cup_{j=1}^{l}X_j= V(G)\\
                        |X_j|=n_j\forall j}
        }
        \prod_{i = 1}^{l}
\left(
q_i!S(m_i,q_i)\ncra{G_{X_i}}{\mathcal{G}(n_i,1,(n_i,q_i))}
\right)\\
 & =    \sum_{\substack{(q_1,q_2,\ldots,q_l)\mid \\
                        n_j-1\leq q_j \leq m_j \forall j}
        }
\left(
        \prod_{i = 1}^{l} q_i!S(m_i,q_i)
\right)
\left(
        \sum_{\substack{(X_1,X_2,\ldots,X_l)\mid \\
                        \cup_{j=1}^{l}X_j= V(G)\\
                        |X_j|=n_j\forall j}
        }
        \prod_{i = 1}^{l}
\ncra{G_{X_i}}{\mathcal{G}(n_i,1,(n_i,q_i))}
\right)
 \\
\end{split}
\end{equation}
The sequence $(n_i,q_i)_{i=1}^l$ may be written as
$(n_i^\prime,q_i^\prime)^{\mu_i}$; $i = 1$ to $r$,
which denotes that the pair $(n_i^\prime,q_i^\prime)$
appears $\mu_i$ times in the sequence $(n_i,q_i)_{i=1}^l$,
the pairs $(n_i^\prime,q_i^\prime)$ are all distinct
for $i = 1$ to $r$, and that $\sum_{i=1}^r \mu_i = l$.
Then each subgraph of $G$ that belongs to the family
$\mathcal{G}(v(G),l,(n_i,q_i)_{i=1}^l)$ is counted
$\prod_{i=1}^r \mu_i!$ times in the inner summation.
Therefore, 
\begin{equation}
\label{eq-dics-subgraphs}
\begin{split}
 T_{v(G)}(A,B\rightarrow G) & =
        \sum_{\substack{(q_1,q_2,\ldots,q_l)\mid \\
                        n_j-1\leq q_j \leq m_j \forall j}
        }
        \left(\prod_{i = 1}^{l}q_i!S(m_i,q_i)\right)
        \left(\prod_{i=1}^r \mu_i!\right)
\ncra{G}{\mathcal{G}(v(G),l,(n_i,q_i)_{i=1}^l)}
\end{split}
\end{equation}
The LHS of Equation~(\ref{eq-dics-subgraphs}) is
known by Theorem~\ref{thm-everything}. 
Now we prove the claim by induction on $\sum_i m_i$.
The base case of induction corresponds to the case
in which each component in the subgraphs being counted
has minimum number of edges, that is,
$m_i = n_i -1$ for all $i\leq l$. In this case,
there is only one term on the RHS of
Equation~(\ref{eq-dics-subgraphs}), and it contains
the unknown $\ncra{G}{\mathcal{G}(v(G),l,(n_i,n_i-1)_{i=1}^l)}$, which can
be solved for. Suppose the claim is true for $\sum_i m_i < m$.
Now let $\sum_i m_i = m$. In this case, as in Lemma~\ref{lem-kedge},
there is only one unknown $\ncra{G}{\mathcal{G}(v(G),l,(n_i,m_i)_{i=1}^l)}$
on the RHS of Equation~(\ref{eq-dics-subgraphs}). All other
terms on the RHS are known by the induction hypothesis.
We can compute the unknown term to obtain the
desired result. This completes the induction step
and the proof.
\end{proof}

\begin{thm}
\label{thm-rank}
The rank polynomial $R(G;x,y)$ is reconstructible from $\mathcal{N}(G)$.
\end{thm}
\begin{proof} Lemma~\ref{lem-lcompo} can be
applied to all induced subgraphs of $G$. So, we
can count the number of subgraphs with $v$ vertices
(none of which isolated), $e$ edges and $l$
components for all $v\leq v(G)$, $e\leq e(G)$ and
$l\geq 1$. Therefore, $\rho_{rs}$ in the expression for the
rank polynomial are known.
\end{proof}


\section{Computing $P(G;\lambda)$ from $\mathcal{PD}(G)$}
\label{sec-poly}
In this section, we consider the problem of computing
the characteristic polynomial of a graph from its
complete polynomial deck. We prove that elementary spanning
subgraphs of each type other than hamiltonian cycles can
be counted from the complete polynomial deck of a graph,
thus proving that the characteristic polynomial of
a non-hamiltonian graph is reconstructible from its
complete polynomial deck.

Here we apply the idea of Kocay's Lemma to the coefficients
of the characteristic polynomials.

Let $A \equiv (a_i)_{i = 1}^k$ be a non-increasing sequence
in $[2,v(G)]$. In this section, we define the notation
$p(A \rightarrow G_X)$ and $c(A\rightarrow G_X)$ differently.

\begin{equation}
\label{eq-pc}
\begin{split}
p(A \rightarrow G_X) & = \prod_{j=1}^{k}(-1)^{a_j}c_{a_j}(G_X)\\  
 & = \prod_{j = 1}^{k}\left(
        \sum_{F\subseteq G_X,\,F\in L_{a_j}}(-1)^{r(F)}2^{s(F)} \right)\\
 & = \sum_{\substack{(F_1,F_2,\ldots,F_{k})\mid \\
               (\forall j\in [1,k])(F_j\subseteq G_X,\,F_j\in L_{a_j})}
            }
        \left((-1)^{\sum_{j=1}^{k} r(F_j)} 2^{\sum_{j=1}^{k} s(F_j)}\right)\\
& = \sum_{Y\subseteq X}c(A\rightarrow G_Y)
\end{split}
\end{equation}
where
\begin{equation}
\begin{split}
c(A\rightarrow G_Y) & =
     \sum_{\substack{(F_1,F_2,\ldots,F_{k})\mid \\
             (\forall j\in [1,k])(F_j\subseteq G_Y,\,F_j\in L_{a_j})\\
                \bigcup_{j=1}^{k}(VF_j)\, =\, Y}
          }
        \left((-1)^{\sum_{j=1}^{k} r(F_j)} 2^{\sum_{j=1}^{k} s(F_j)}\right)
\end{split}
\end{equation}
Thus we have grouped together tuples of elementary 
subgraphs in groups that span each subset of $X$.
\begin{lem}
\label{lem-pd2cdec} If $A$ is a sequence defined over $[2,v(G)-1]$
then $c(A\rightarrow G)$ is reconstructible from the complete polynomial
deck of $G$.
\end{lem}

\begin{proof} As in the proof of Lemma~\ref{lem-exp-c},
$p(A \rightarrow G_X)$ can be computed for each induced subgraph
$G_X$. By M\"{o}bius inversion of Equation~(\ref{eq-pc}),
we write
\begin{equation}
\label{eq-pd-mobius}
c(A\rightarrow G_X) = \sum_{Y\subseteq X}(-1)^{|X\backslash Y|}
p(A\rightarrow G_Y)
\end{equation}
But we cannot compute the RHS of Equation~(\ref{eq-pd-mobius})
because, in general, for $X\subsetneq VG$,
we do not know which polynomials in $\mathcal{PD}(G)$
correspond to the induced subgraphs of $G_X$. But this is
not a problem if $X = VG$.
We can write
\begin{equation} 
c(A\rightarrow G) = \sum_{Y\subseteq VG}(-1)^{|VG\backslash Y|}
p(A\rightarrow G_Y)
\end{equation}
Now the RHS, and hence $c(A\rightarrow G)$, can be computed.
\end{proof}

\noindent {\bf Remark.} Note that we would not be
able to compute $p(A\rightarrow G)$ if we defined
$A$ in $[2,v(G)]$, because we do not know $c_{v(G)}(G)$.

Let $\lambda (m,p) \equiv (x_1, x_2, \ldots , x_p)$ denote a
partition of $m$. We assume that $x_1\, \geq \, x_2\, \ldots \,
\geq x_p$. We write $\lambda (m) $ when the number of parts $p$
is not relevant. Also, we just write $\lambda $ instead of 
$\lambda (m,p)$ when $m$ and $p$ are either understood from the context
or not relevant. Another partition $\lambda ^\prime (m,p+1) $
can be obtained from $\lambda (m,p)$ by replacing an
$x_i$ by $y$ and $z$ such that $y+z = x_i$, and ordering
the numbers in a non-increasing order. Any partition that is
obtained from $\lambda (m,p)$ by a sequence of such operations
is called a {\em refinement} of $\lambda (m,p)$. Also,
$\lambda (m,p)$ is a trivial refinement of itself.
If $\lambda ^\prime (m,q)$ is a refinement of $\lambda (m,p)$,
then we denote it by $\lambda ^\prime (m,q) \preceq \lambda (m,p)$.
The relation $\preceq $ between partitions is a partial order.

Now consider partitions in which the smallest part
$x_p$ is at least 2. Associated with each such partition
$\lambda (m,p)$, there is a unique elementary graph
$F_{\lambda }\in L_m$, whose $i$'th component is a
cycle of length $x_i$, or an edge if $x_i\, =\, 2$.
If the sequence $A\equiv (a_i)_{i=1}^k$ is defined
such that $\lambda \equiv (a_1,a_2, \ldots, a_k)$ is
a non-trivial partition of $v(G)$, then we denote
$c(A\rightarrow G)$ by $c(\lambda \rightarrow G)$.

\begin{lem}
\label{lem-refinement}
Let $\lambda  \equiv (a_1,a_2, \ldots, a_k)$
be a non-trivial partition of $n = v(G)$. Then,
\begin{equation}
\label{eq:refinement}
\begin{split}
c(\lambda \rightarrow G) & =
        \sum_{\lambda ^\prime \preceq \lambda}
        c(\lambda \rightarrow F_{\lambda ^\prime})
        \ncra{G}{F_{\lambda ^\prime}}
\end{split}
\end{equation}
\end{lem}
\begin{proof} From the definition of $c(A\rightarrow G)$,
we write
\begin{equation}
\begin{split}
c(\lambda\rightarrow G) & =
\sum_{\substack{(F_1,F_2,\ldots,F_{k})\mid \\
        (\forall j\in [1,k])(F_j\subseteq G,\,F_j\in L_{a_j})\\
                \bigcup_{j=1}^{k}(VF_j)\, =\, VG}
     }
        \left((-1)^{\sum_{j=1}^{k} r(F_j)} 2^{\sum_{j=1}^{k} s(F_j)}\right)\\
& = \sum_{F\subseteq G, F\in L_n}
\sum_{\substack{(F_1,F_2,\ldots,F_{k})\mid \\
        (\forall j\in [1,k])(F_j\subseteq F,\,F_j\in L_{a_j})\\
                 \bigcup_{j=1}^{k}F_j\, =\, F}
}
\left((-1)^{\sum_{j=1}^{k} r(F_j)} 2^{\sum_{j=1}^{k} s(F_j)}\right)\\
& =  \sum_{\substack{F\subseteq G,F\in L_n}} c(\lambda \rightarrow F)\\
& = \sum_{\lambda ^\prime \preceq \lambda}
        c(\lambda \rightarrow F_{\lambda ^\prime})
        \ncra{G}{F_{\lambda ^\prime}}
\end{split}
\end{equation}
The last step may be explained as follows:
if $F$ is a disjoint union of elementary graphs
$F_j\in L_{a_j}; j \in [1,k]$, where $v(F) = v(G) = n$,
then $F$ is isomorphic to an elementary graph $F_{\lambda ^\prime}$
for some refinement ${\lambda ^\prime}$ of $\lambda$. Trivially,
if each $F_j$ is the cycle $C_{a_j}$ then $F = F_{\lambda}$.
We then group the terms $c(\lambda \rightarrow F)$ by the
isomorphism type of $F$.
\end{proof}

\begin{lem}
\label{lem-non-hamiltonian}
If $F$ is an elementary graph on $n = v(G)$ vertices,
other than the cycle, then $\ncra{G}{F}$
is reconstructible from $\mathcal{PD}(G)$. 
\end{lem}

\begin{proof}
Since $F$ is not a cycle, it is isomorphic to
$F_{\lambda_0}$ for a unique non-trivial
partition $\lambda_0 $ of $v(G)$.
From Equation~(\ref{eq:refinement}) we write
\begin{equation}
\label{eq:nonh1}
\begin{split}
\ncra{G}{F_{\lambda_0}} &= \frac{1}{c(\lambda_0 \rightarrow F_{\lambda_0 })}
\left(
c(\lambda_0 \rightarrow G)
     -\sum_{\lambda \prec \lambda_0}
        c(\lambda_0 \rightarrow F_{\lambda })
        \ncra{G}{F_{\lambda}}
\right)
\end{split}
\end{equation}
Now we expand $ \ncra{G}{F_{\lambda}} $ on the RHS
of the above equation by repeated application of
the same equation, and obtain the following solution.
\begin{equation}
\label{eq:nonh2}
\ncra{G}{F_{\lambda_0}}
 = \sum_{\lambda_q \prec \lambda_{q-1} \prec \ldots \prec \lambda_0}
\frac{\left(-1\right)^q c(\lambda_q \rightarrow G)
\prod_{i=0}^{q-1}c(\lambda_i\rightarrow F_{\lambda_{i+1}})}
{\prod_{i=0}^qc(\lambda_i\rightarrow F_{\lambda_i})}
\end{equation}
where the summation is over all chains
$\lambda_q \prec \lambda_{q-1} \prec \ldots \prec \lambda_0$; 
$q \geq 0$, and an empty product is 1. There are finitely
many terms in the above summation since there are
finitely many refinements of $\lambda_0$. Since $\lambda_0 $ is a
non-trivial partition of $n$, (that is, $x_i < v(G)\,\forall\, i$),
by Lemma~\ref{lem-pd2cdec}, $c(\lambda \rightarrow G)$
is reconstructible for each $\lambda \preceq \lambda_0$.
Also, for each $\lambda \preceq \lambda_0$,
$c(\lambda \rightarrow F_{\lambda })$ is non-zero.
(Here we would like to repeat that we have considered only those
partitions in which the smallest part is at least 2.)
Thus the RHS can be computed.
\end{proof}

The main theorem in this section now follows from the above lemmas.
\begin{thm}
\label{thm-deg1}
If $F$ is an elementary graph on $v(G)$ vertices,
other than a cycle, then $\ncra{G}{F}$ is
reconstructible from $\mathcal{PD}(G)$.
Therefore, if there is a vertex of degree 1 in $G$,
then the characteristic polynomial of $G$ is reconstructible
from its complete polynomial deck.
\end{thm}

\begin{proof}
The degree sequence of a graph is reconstructed
from its complete polynomial deck as follows.
Consider the polynomials of degree $v(G)-1$ in
$\mathcal{PD}(G)$. They are the characteristic polynomials of
the vertex deleted subgraphs $G-u$ of $G$ for $u\in VG$.
Since $c_2(G)$ and  $c_2(G-u)$ count the number of edges
of $G$ and $G-u$, respectively, we know the degree
of $u$ in $G$ for each $u \in VG$.  
Thus the premise of Theorem~\ref{thm-deg1} is recognised from
$\mathcal{PD}(G)$. 
The coefficients $c_i(G);i<v(G)$ can be computed
using Lemma~\ref{lem-derivative}. 
Since there are no hamiltonian cycles in $G$,
Lemma~\ref{lem-non-hamiltonian} implies that the constant term in
the characteristic polynomial of $G$ can be calculated.
\end{proof}

\noindent {\bf Remark.} Whenever non-hamiltonicity of a graph
is recognised from its complete polynomial deck, its
characteristic polynomial can be computed as well.


\section{Whitney's Theorem}
\label{sec-whitney}
In Section~\ref{sec-intro}, it was stated that the
computation of the chromatic polynomial of a graph requires only
non-separable induced subgraphs of the graph.
Whitney's proof of this fact was
based on his theorem that separable subgraphs can
be counted from the counts of non-separable subgraphs. 
Let $n_t(G)$ denote the number of subgraphs of $G$
of type $t$, where `type' of a graph is determined
by the number of blocks of each isomorphism type.
He proved, (stated in the terminology of \cite{biggs1993}),
that there is a polynomial $\phi_t$, independent of $G$, such that
\begin{equation}
n_t(G) = \phi_t(n_{\sigma}(G), n_{\rho}(G), ...)
\end{equation}
where $\sigma$, $\rho$, ... are non-separable types
with not more edges than $t$.
Here we prove \mbox{Whitney's} result using Kocay's Lemma.
Our presentation explicitly describes the poly\-nomial in Whitney's
theorem.

Let $S_0 = \{F_1, F_2, \ldots ,F_k\}$ be a family of non-separable
graphs, some of them possibly isomorphic. Thus $S_0$ represents
a `graph type'. Extending the notation introduced earlier,
we write $\ncra{G}{S_0}$ to denote the number of subgraphs
of $G$ of type $S_0$.
We define a partial order $\preceq $ on graph types as follows.
Let $S$ be a graph type. We say that $S\preceq S_0$ if $c(S_0,X)$
is non-zero for some graph $X$ of type $S$. It is easily
seen that $c(S_0,X)$ depends only on the type $S$ of $X$,
not on a particular choice of $X$. So, we write it as $c(S_0,S)$.
For any graph $G$, by Kocay's Lemma~\ref{lem-kocay},
\begin{equation}
p(S_0) = \prod_{i = 1}^{k} \ncra{G}{F_i} = \sum_{X} c(S_0,X)\ncra{G}{X}
\end{equation}
Terms on the RHS can be grouped together according to
the types of $X$, so we can write
\begin{equation}
\begin{split}
p(S_0) &= \sum_{S\preceq S_0} c(S_0,S)\ncra{G}{S} \\
&= c(S_0,S_0)\ncra{G}{S_0}+\sum_{S\prec S_0} c(S_0,S)\ncra{G}{S}
\end{split}
\end{equation}
Therefore,
\begin{equation}
\ncra{G}{S_0} = \frac{1}{c(S_0,S_0)}
\left(
p(S_0)-\sum_{S\prec S_0}c(S_0,S)\ncra{G}{S}
\right)
\end{equation}
Now we repeatedly apply the same equation to $\ncra{G}{S}$
on the RHS, as we did in Equation~(\ref{eq:nonh2}). 
We thus get the polynomial of Whitney's theorem.
\begin{thm}
\label{thm-whitney1}
\begin{equation}
\label{eq-whitney}
\ncra{G}{S_0} = \sum_{S_q \prec S_{q-1} \prec \ldots \prec S_0}
\frac{\left(-1\right)^qp(S_q)\prod_{i=0}^{q-1}c(S_i,S_{i+1})}
{\prod_{i=0}^qc(S_i,S_i)}
\end{equation}
where the summation is over all chains
$S_q \prec S_{q-1} \prec \ldots \prec S_0$; $q \geq 0$,
and an empty product is 1.
\end{thm}
There are finitely many terms in the summation in
Equation~(\ref{eq-whitney}) because each $S_q$ has fewer
blocks in it than $S_{q-1}$.
While some other known proofs of
this theorem are based on not very different ideas,
(for example, see \cite{biggs1978}),
the above explicit formulation of the polynomial seems new.
It allows us to argue about the reconstructibility of
the characteristic polynomial more directly than in other
standard proofs.

\begin{cor}
\label{cor-w2p}
The characteristic polynomial of a graph is reconstructible from
its vertex deck.
\end{cor}

\begin{proof} Let $G$ be the graph under consideration.
Its elementary spanning subgraphs other
than the hamiltonian cycles are counted as
in the standard proof.
Let $n = v(G)$, and let $S_0$ be the type of an
$n$-vertex elementary graph $H$ other than $C_n$.
Any block in a type $S \prec S_0$ has
fewer than $n$ vertices, and $H$ is the only graph
of type $S_0$ that has $n$ vertices. So $\ncra{G}{H}$ can be
counted using Whitney's Theorem~\ref{thm-whitney1}
and Kelly's Lemma~\ref{lem-kelly}.

To count hamiltonian cycles, we set $S_0 = \{nK_2\}$, that
is, a graph type consisting of $n$ blocks, each one of them
a $K_2$. Since no subgraph of $G$ has $n$ blocks isomorphic
to $K_2$, $\ncra{G}{S_0} = 0$. On the RHS of Equation
~(\ref{eq-whitney}), there is precisely one term
that contributes hamiltonian cycles. That is,
$S_1 = \{C_n\} \prec S_0$ is the unique chain
that contributes $ham(G)$, implying that
the terms containing $ham(G)$ cannot cancel out.
All other blocks that appear in Equation~(\ref{eq-whitney})
have fewer than $n$ vertices,
so can be counted using Kelly's Lemma~\ref{lem-kelly}.
Therefore, $ham(G)$ is reconstructible.
The reconstructibility of $P(G;\lambda)$ now follows from
Lemma~\ref{lem-sachs}.
\end{proof}

\noindent {\bf Remark.} In the standard proof of the
reconstructibility of the characteristic polynomial,
one applies Kocay's Lemma directly.
As a result one has to proceed step by step, counting
spanning trees, then spanning unicyclic subgraphs, etc.,
as we did in Corollary~\ref{cor-trees}. These intermediate
steps are skipped by the direct application of Whitney's
theorem.

\section{Problems and discussion}
Expressing $c_n(G)$ or $\psi_n(G)$, where $n = v(G)$,
as polynomials in $c_j(G-S)$ or $\psi_j(G-S)$;
$S\subsetneq VG$, would be of interest.
Alternatively, we would like to construct a generalisation
of the characteristic polynomial which can be computed more naturally
from the poset of induced subgraphs, and from which the
characteristic polynomial can be easily computed.
Such a goal is motivated on the one hand by the proofs in
Section~\ref{sec-nmatrix}, and, on the other hand,
by similar generalisations of the chromatic polynomial,
viz, Stanley's chromatic symmetric function, (see
\cite{stanley1995} \& \cite{stanley1999}),
and another recent two variable generalisation of the chromatic
polynomial \cite{dohmen2003}. Both these generalisations
are closely related to the lattice of connected partitions
of $VG$, (see \cite{stanley1995} for definitions). A
relationship between the poset of induced subgraphs defined
in this paper and the lattice of connected partitions of
the vertex set defined by Stanley could possibly be established
using Kocay's Lemma. Such a result would be
useful in understanding exact relationship between
different expansions (and reconstructibility)
of several important invariants in a unified way.

The reconstruction of the number of hamiltonian cycles is
difficult and indirect in the proofs we have
presented here, and in the original proof by Tutte as well.
A reason for this difficulty is seen in Whitney's theorem.
Observe that for an $n$-vertex graph $G$, the polynomial
of Whitney's theorem contains a term in $ham(G)$ only
if the type $S_0$ contains $n$ copies of $K_2$, or a $C_n$,
and possibly other types of blocks. As a result,
in Corollary~\ref{cor-w2p} we had to count all possible
blocks with at most $n$ edges. But can we count the
number of hamiltonian cycles from the $\mathcal{N}$-matrix
at least as clearly as in Corollary~\ref{cor-w2p}? Towards
this goal, we would like to understand the
relationship between the structure of the edge labelled
poset for separable graphs and that for blocks,
and prove a generalisation of Whitney's theorem.

The crucial difference between the proofs in
Section~\ref{sec-nmatrix} and Section~\ref{sec-poly}
is in Lemmas~\ref{lem-exp-c} and ~\ref{lem-pd2cdec}.
In Lemma~\ref{lem-pd2cdec}, the use of M\"{o}bius inversion
was limited to the computation of $c(A\rightarrow G)$ because
we did not know the partial order on the induced subgraphs.
This suggests that a general `expansion' for the number
of hamiltonian cycles would probably involve a summation over
chains in $\mathcal{ELP}(G)$. Therefore, counting hamiltonian
cycles from $\mathcal{PD}(G)$ and the original problem of
Gutman and Chvetkovi\'{c} seem difficult. This is
probably why many known results on the reconstruction
of the characteristic polynomial of a graph from its
characteristic polynomial deck assume the graph to contain
several pendant vertices.

We propose the following generalisation of reconstruction
for studying questions simi\-lar to the one posed by
Gutman and Chvetkovi\'{c}. 
Suppose $f$ is a graph invariant, and we are interested
in reconstructing $G$ or partial invariants of $G$ from the
deck $\mathcal{D}(G;f) = \{f(G-u); u \in VG\}$. A new
collection $\mathcal{D}!(G;f)$ is recursively defined
as $\{(f(G-u), \mathcal{D}!(G;f-u)); u\in VG\}$.
We then define an equivalence relation $\sim $
on graphs such that $H_1 \sim H_2$ if
$(f(H_1), \mathcal{D}!(H_1)) = (f(H_2),\mathcal{D}!(H_2))$.
This relation gives an incidence matrix (or an edge labelled
poset) on the types of induced subgraphs
of $G$, where `type' refers to an equivalence
class under the relation $\sim $ defined above.
It can be shown that for many invariants $f$,
Ulam's conjecture is true if and only if all graphs $G$ on
more than 2 vertices are reconstructible from
$\mathcal{D}!(G;f)$. One example of such an invariant is:
$f(G) = 1$ if $G$ has a vertex of degree 1, and $f(G) = 0$
otherwise. Another example is $f(G) = P(G;\lambda)$.
The proof of this is similar to that of
Proposition~\ref{prop:equiv-un}: the base case follows from the fact
that any three vertex graph $G$ is completely determined
by $\mathcal{D}!(G;f)$ for the above invariants.
The problem of reconstructing $G$ from the deck $\mathcal{D}!(G;f)$
is similar to the generalisation of the reconstruction
problem suggested by Tutte, (Notes on pp. 123-124 in
\cite{tutte1984}). We are not really interested in
the question of computing $f(G)$ from $\mathcal{D}(G;f)$.
Rather we ask the question - what are the incomplete
invariants $f$, (that is, the invariants that do not
determine a graph completely,) and classes of graphs $G$,
for which $\mathcal{D}!(G;f)$ could be constructed from
$\mathcal{D}(G;f)$? If we could construct $\mathcal{D}!(G;f)$
from $\mathcal{D}(G;f)$, then we could also prove all
the results of Section~\ref{sec-nmatrix}. We would like
to investigate this question when $f(G) = P(G;\lambda)$,
and when $f(G) = (P(G;\lambda), P(\bar{G};\lambda))$
- the invariant which was considered by Hagos 
\cite{hagos2000}.

\subsection*{Acknowledgements} I take this opportunity
to thank Allan Wilson Centre for the support and
encouragement. I would also like to thank the referee
for several useful suggestions for improving the
presentation.

\end{document}